\documentclass[12pt]{article}
\usepackage[utf8]{inputenc}
\usepackage{amsmath, amssymb, url, float, graphicx, float}
\usepackage{fullpage}
\usepackage[small,bf]{caption}
\usepackage{cite}
\usepackage{comment}
\usepackage[colorlinks=true, citecolor=blue, linkcolor=blue]{hyperref}
\usepackage{xcolor}
\usepackage{listings}

\definecolor{mediumviolet-red}{rgb}{0.78, 0.08, 0.52}

\lstdefinelanguage{mypython}
{
	keywords=[1]{from, import, assert, not, print},
	keywordstyle=[1]{\color{mediumviolet-red}},
	keywords=[2]{pymde, torch},
	keywordstyle=[2]{\color{blue}},
	numbers=none,
	upquote=true,
	showstringspaces=false,
	basicstyle=\ttfamily,
	columns=fullflexible,
	keepspaces=true,
	emph={True,False,as,def,return,float},
	emphstyle={\color{blue}},
	frame=trbl,
	belowskip=1em,
	aboveskip=1em,
	captionpos=b
}

\renewcommand\phi\varphi

\let\phi\varphi

\newcommand{\ones}{\mathbf 1}
\newcommand{\reals}{{\mbox{\bf R}}}

\newcommand{\Expect}{\mathop{\bf E{}}}

\newcommand{\eg}{{\it e.g.}}
\newcommand{\ie}{{\it i.e.}}

\newcommand{\BEAS}{\begin{eqnarray*}}
\newcommand{\EEAS}{\end{eqnarray*}}
\newcommand{\BEA}{\begin{eqnarray}}
\newcommand{\EEA}{\end{eqnarray}}
\newcommand{\BEQ}{\begin{equation}}
\newcommand{\EEQ}{\end{equation}}
\newcommand{\BIT}{\begin{itemize}}
\newcommand{\EIT}{\end{itemize}}

\title{Constant Function Market Makers: Multi-Asset Trades via Convex Optimization}
\author{Guillermo Angeris\\\texttt{\small angeris@stanford.edu} \and
Akshay Agrawal \\\texttt{\small akshayka@stanford.edu}\and
Alex Evans\\\texttt{\small ahe4nc@gmail.com} \and
Tarun Chitra \\ \texttt{\small tarun@gauntlet.network} \and
Stephen Boyd \\ \texttt{\small boyd@stanford.edu}}
\date{July 2021}

\begin{document} 
\maketitle 

\begin{abstract}
The rise of Ethereum and other blockchains that support smart contracts has led
to the creation of decentralized exchanges (DEXs), such as
Uniswap, Balancer, Curve, mStable, and SushiSwap, which enable agents to trade
cryptocurrencies without trusting a centralized authority.
While traditional exchanges use order books to match and execute trades,
DEXs are typically organized as constant function market makers (CFMMs).
CFMMs accept and reject proposed trades based on the evaluation of a 
function that depends on the proposed trade and the current reserves of the DEX.
For trades that involve only two assets, CFMMs are easy to understand,
via two functions that give the quantity of one asset that must be tendered to
receive a given quantity of the other, and vice versa.
When more than two assets are being exchanged, it is harder to
understand the landscape of possible trades.  
We observe that various problems of choosing
a multi-asset trade can be formulated as convex optimization problems,
and can therefore be reliably and efficiently solved.
\end{abstract}

\section{Introduction}
In the past few years, several new financial exchanges have been implemented
on blockchains, which are distributed and permissionless ledgers replicated
across networks of computers.
These \emph{decentralized exchanges} (DEXs) enable agents to trade
cryptocurrencies, \ie, digital currencies with account balances stored on a
blockchain, without relying on a trusted third party to facilitate the 
exchange. DEXs have significant capital flowing through them; the
four largest DEXs on the Ethereum blockchain (Curve Finance
\cite{egorovStableSwapEfficientMechanism}, Uniswap \cite{uniswap,
adams2021uniswap}, SushiSwap \cite{sushiswap}, and Balancer \cite{balancer})
have a collective trading volume of several billion dollars per day.

Unlike traditional exchanges, DEXs typically do not use
order books. Instead, most DEXs (including Curve, Uniswap, SushiSwap, and
Balancer) are organized as \emph{constant function market makers} (CFMMs). A
CFMM holds reserves of assets (cryptocurrencies), contributed by liquidity
providers. Agents can offer or tender baskets of assets to the CFMM, in
exchange for another basket of assets. If the trade is accepted, the tendered
basket is added to the reserves, while the basket received by the agent is subtracted
from the reserves. Each accepted trade incurs a small fee, which is distributed
pro-rata among the liquidity providers.

CFMMs use a single rule that determines whether or not a proposed trade is
accepted. The rule is based on evaluating a \emph{trading function}, which
depends on the proposed trade and the current reserves of the CFMM. A proposed
trade is accepted if the value of the trading function at the post-trade
reserves (with a small correction for the trading fee)
equals the value at the current reserves, \ie, the function is held
constant. This condition is what gives CFMMs their name. One simple example of
a trading function is the product \cite{lu2017building, buterin2017path},
implemented by Uniswap \cite{uniswap} and SushiSwap
\cite{sushiswap}; this CFMM accepts a trade only if it leaves the product of
the reserves unchanged. Several other functions can be used, such as the sum or
the geometric mean (which is used by Balancer \cite{balancer}). 
%As we will see 
%later, the choice of function depends on many things, such as the types of
%assets being traded and their expected volatility.

For trades involving just two assets, CFMMs are very simple to understand,
via a scalar function that relates how much of one asset is required to receive 
an amount of the other, and vice versa.  Thus the choice of a two-asset
trade involves only one scalar quantity: 
how much you propose to tender (or, equivalently,
how much you propose to receive).

For general trades, in which many assets may be simultaneously exchanged,
CFMMs are more difficult reason about. 
When multiple assets are tendered,
there can be many baskets that can be tendered to receive
a specific basket of assets, and vice versa, there are many choices
of the received basket, given a fixed one that is tendered.
Thus the choice of a multi-asset trade is more complex than just
specifying an amount to tender or receive.
In this case the trader may wish to tender and receive baskets that are
most aligned with their preferences or 
utility (\eg, one that maximizes their risk-adjusted return).

In all practical cases, including the ones mentioned above, the trading
function is concave~\cite{angerisImprovedPriceOracles2020}. In this paper we
make use of this fact to formulate various multi-asset
trading problems as convex optimization problems. Because convex optimization
problems can be solved reliably and efficiently (in theory and in practice)
\cite{boyd2004convex}, we can solve the formulated trading problems exactly.
This gives a practical solution to the problem of choosing among many
possible multi-asset trades: the trader articulates their
objective and constraints, and a solution to this problem determines the baskets
of assets to be tendered and received.

\paragraph{Outline.} 
We start by surveying related work in \S\ref{s-rel-work}. In \S\ref{s-cfmm}, we
give a complete description of CFMMs, describing how agents may trade with a
CFMM, as well as add or remove liquidity. In \S\ref{s-prop} we study some basic
properties of CFMMs, many of which rely on the concavity of the trading function.
In \S\ref{s-two-asset} we examine trades involving just two assets,
and show how to understand them via two functions that give the amount of asset
received for a given quantity of the tendered asset. Finally, in
\S\ref{s-multi-asset} we formulate the general multi-asset trading problem as a
convex optimization problem, and give some specific examples.

\subsection{Background and related work}\label{s-rel-work}

\paragraph{Blockchain.} CFMMs are typically implemented on a \emph{blockchain}:
a decentralized, permissionless, and public ledger. The blockchain stores 
accounts, represented by cryptographic public keys, and associated balances of
one or more cryptocurrencies. A blockchain allows any two accounts to securely
transact with each other without the need for a trusted third party or central
institution, using public-key cryptography to verify their identities. Executing
a \emph{transaction}, which alters the state of the blockchain, costs the
issuer a fee, typically paid out to the individuals providing computational
power to the network.
(This network fee depends on the amount of computation a transaction requires and is paid 
in addition to the CFMM trading fee mentioned above and described below.)

Blockchains are highly tamper resistant: they
are replicated across a network of computers, and kept in consensus via simple
protocols that prevent invalid transactions such as double-spending of a coin.
The consensus protocol operates on the level of \emph{blocks} (bundles of
transactions), which are verified by the network and chained together to form
the ledger. Because the ledger is
public, anyone in the world can view and verify all account balances and the
entire record of transactions. 

The idea of a blockchain originated with a pseudonymously authored
whitepaper that proposed Bitcoin, widely considered to be the first
cryptocurrency~\cite{nakamoto2008bitcoin}. 

\paragraph{Cryptocurrencies.} A cryptocurrency is a digital currency
implemented on a blockchain. Every blockchain has its own native
cryptocurrency, which is used to pay the network transaction fees (and can also
be used as a standalone currency).

A given blockchain may have several other cryptocurrencies implemented on it.
These additional currencies are sometimes called \emph{tokens}, to distinguish
them from the base currency. There are thousands of tokens in circulation
today, across various blockchains. Some, like the Uniswap token UNI, give
holders rights over the governance of a protocol, while others, like USDC, are
\emph{stablecoins}, pegged to the market value of some external or real-world
currency or commodity.

\paragraph{Smart contracts.} 
Modern blockchains, such as Ethereum
\cite{buterin2013ethereum, wood2014ethereum}, Polkadot \cite{wood2016polkadot},
and Solana \cite{yakovenko2018solana}, allow anyone to deploy
arbitrary stateful programs called \emph{smart contracts}. A contract's public
functions can be invoked by anyone, via a transaction sent through the network
and addressed to the contract. (The term
`smart contract' was coined in the 1990s, to refer to a set of 
promises between agents codified in a computer program \cite{szabo1995smart}.)
%A smart contract is
%a special type of account. Like regular accounts, they have an
%associated balance, but unlike regular accounts, they carry arbitrary programs
%(code and state) with publicly accessible functions.
%Blockchains that support smart contracts can be viewed as distributed global
%computers, with open source code and open state.
Because creators are free to compose deployed contracts or remix them in their
own applications, software ecosystems on these blockchains have developed
rapidly.

CFMMs are implemented using smart contracts, with functions for trading, adding
liquidity, and removing liquidity. Their implementations are usually simple.
For example, Uniswap v2 is implemented in just 200 lines of code. In addition
to DEXs, many other financial applications have been deployed on blockchains,
including lending protocols (\eg, \cite{aave, compound}) and various derivatives
(\eg, \cite{uma, dydx}). The collection of financial applications running on
blockchains is known as decentralized finance, or DeFi for short.

\paragraph{Exchange-traded funds.} 
CFMMs have some similarities to exchange-traded funds (ETFs). A CFMM's
liquidity providers are analogous to an ETF's authorized participants;
adding liquidity to a CFMM is analogous to the creation of an ETF share,
and subsequently removing liquidity is analogous to redemption. But while the
list of authorized participants for an ETF is typically very small, anyone in
the world can provide liquidity to a CFMM or trade with it.

%\paragraph{Comparison to order books.} XXX I don't understand any of this, which
%has nothing to with anything in this paper, as far as I can tell. XXX
%One of the main advantages of CFMMs over order books is composablity. One can
%compose a trade from asset A to asset B with another trade from asset B to
%asset C to generate a trade from asset A to asset C. Unlike order books, CFMMs
%enable compositions of trades to be computed in constant time (relative to the
%number of liquidity providers). In compute-constrained blockchain environments,
%CFMMs have had orders of magnitude more liquidity and trading volume relative
%to order book models (e.g. Etherdelta, Serum) in part due to their low compute
%and storage requirements. When coupled with models for atomic trade execution
%(e.g. flash loans), CFMMs make it easy execute a trade from A to B to C to D.
%The routing of trades has allowed for millions of illiquid assets to be
%tradeable with reasonably low trading costs relative to order books.

\paragraph{Comparison to order books.}
In an order book, trading a basket of multiple assets for another basket of 
multiple assets requires multiple separate trades.  Each of these
trades would entail the blockchain fee, increasing the total cost of 
trading to the trader.  In addition, multiple trades
cannot be done at the same time with an order book,
exposing the trader to the risk that some of 
the trades go through while others do not, or that some of the trades will execute 
at unfavorable prices.
In a CFMM, multiple asset baskets are exchanged in one trade, which either 
goes through as one group trade, or not at all, so the trader is not exposed 
to the risk of partial execution.

Another advantage of CFMMs over order book exchanges is their efficiency of 
storage, since they do not
need to store and maintain a limit order book, and their computational efficiency,
since they only need to evaluate the trading function. Because users must pay for
computation costs for each transaction, and these costs can often be nonnegligible in
some blockchains, exchanges implementing CFMMs can often be much cheaper for users to interact with
than those implementing order books.

\paragraph{Previous work.}
Academic work on automated market makers began with the study of scoring rules within the
statistics literature, \eg,~\cite{winkler1969scoring}.
Scoring rules furnish probabilities for baskets of events, which can be
viewed as assets or tokens in a prediction market.
The output probability from a scoring rule was first proposed as a pricing mechanism for a
binary option (such as a prediction market) in~\cite{hanson2003combinatorial}.
Unlike CFMMs, these early automated market makers were shown to be computationally
complicated for users to interact with. For example.
Chen~\cite{chen2008complexity} demonstrated that computing optimal arbitrage portfolios in
logarithmic scoring rules (the most popular class of scoring rules) is \#P-hard.

The first CFMM on Ethereum (the most commonly used blockchain for smart
contracts) was Uniswap~\cite{uniswap, adams2021uniswap}.
The first formal analysis of Uniswap was first done in~\cite{angerisAnalysisUniswapMarkets2020} and
extended to general concave trading functions in~\cite{angerisImprovedPriceOracles2020}.
Evans~\cite{evans2020liquidity} first proved that constant mean market makers could replicate a large set of
portfolio value functions.
The converse result was later proven, providing a mechanism for constructing a trading
function that replicates a given portfolio value function~\cite{angeris2021replicating}.
Analyses of how fees~\cite{evans2021optimal, tassy2020growth} and trading function
curvature~\cite{angeris2020does, aoyagiLiquidityProvisionAutomated, aoyagiLiquidityImplicationsConstant}
affect liquidity provider returns are also common in the literature.
Finally, we note that there exist investigations of privacy in CFMMs~\cite{angeris2021note}, 
suitability of liquidity provider shares as a collateral asset~\cite{chitra2021note},
and the question of triangular arbitrage~\cite{wang2021cyclic} in CFMMs.

\subsection{Convex analysis and optimization} 
\paragraph{Convex analysis.} 
A function $f : D \to \reals$, with $D \subseteq \reals^n$,
is convex if $D$ is a convex set and
\[
    f(\theta x + (1 - \theta)y) \leq \theta f(x) + (1 - \theta)f(y),
\]
for $0 \le \theta \le 1$ and all $x,y \in D$.
It is common to extend a convex function to an extended-valued function
that maps $\reals^n$ to $\reals \cup \{\infty\}$,
with $f(x)=+\infty$ for $x \not\in D$.
A function $f$ is concave if $-f$ is convex \cite[Chap.\ 3]{boyd2004convex}. 

When $f$ is differentiable, an equivalent characterization of convexity is
\[
f(z) \geq f(x)+ \nabla f(x)^T(z-x),
\]
for all $z,x \in D$.
A differentiable function $f$ is concave if and only if
for all $z,x\in D$ we have
\BEQ\label{e-concavity}
f(z) \leq f(x)+ \nabla f(x)^T(z-x).
\EEQ
The right hand side of this inequality is the first-order Taylor 
approximation of the function $f$ at $x$, so this inequality
states that for a concave function, the Taylor approximation
is a global upper bound on the function. 

By adding \eqref{e-concavity} and the same inequality with $x$ and $z$ swapped,
we obtain the inequality 
\BEQ\label{e-monotone}
(\nabla f(z)-\nabla f(x))^T (z-x) \leq 0,
\EEQ
valid for any concave $f$ and $z,x\in D$.
This inequality states that for a concave function $f$, $-\nabla f$ is a 
monotone operator~\cite{ryuPrimerMonotoneOperator2016}.

\paragraph{Convex optimization.}
A \emph{convex optimization problem} has the form
\begin{equation*}
\begin{array}{ll}
\mbox{minimize} & f_0(x) \\
\mbox{subject to} & f_i(x) \leq 0, \quad i=1, \ldots, m\\
& g_i(x) = 0, \quad i=1, \ldots, p,
\end{array}
\end{equation*}
where $x \in \reals^n$ is the optimization variable, the objective function
$f_0 : D \to \reals$ and inequality constraint functions $f_i : D
\to \reals$ are convex, and the equality constraint functions $g_i : \reals^n
\to \reals$ are affine, \ie, have the form $g_i(x) = a_i^T x + b_i$ for some
$a_i \in \reals^n$ and $b_i \in \reals$.
(We assume the domains of the objective and inequality functions
are the same for simplicity.)
The goal is to find a \emph{solution} of the problem, which is
a value of $x$ that minimizes the objective function, among all $x$ satisfying
the constraints $f_i(x) \leq 0$, $i=1, \ldots, m$, and $g_i(x) = 0$, $i=1, \ldots, p$
\cite[Chap.\ 4]{boyd2004convex}. 
In the sequel we will refer to the problem of maximizing a concave function, subject to
convex inequality constraints and affine equality constraints,
as a convex optimization problem, since this problem
is equivalent to minimizing $-f_0$ subject to the constraints.

Convex optimization problems are notable because they have many applications,
in a wide variety of fields,
and because they can be solved reliably and efficiently \cite{boyd2004convex}.
The list of applications of convex optimization is
large and still growing.  It has applications in vehicle control
\cite{stewart2008model, spacex, lipp2014minimum}, finance \cite{cornuejols2006,
boyd2017multi}, dynamic energy
management \cite{moehle2018dynamicnot}, resource allocation
\cite{agrawal2021allocation}, machine learning \cite{friedman2001elements,
boyd2011distributed}, inverse design of physical systems
\cite{angerisHeuristicMethodsPerformance2021}, circuit design
\cite{hershenson2001opamp, boyd2005digital}, and many other fields.

In practice, once a problem is formulated as a convex optimization problem, we
can use off-the-shelf solvers (software implementations of numerical
algorithms) to obtain solutions. Several solvers, such as OSQP
\cite{stellatoOSQPOperatorSplitting2020}, SCS
\cite{odonoghueConicOptimizationOperator2016}, ECOS \cite{ecos}, and COSMO
\cite{garstkaCOSMOConicOperator2019}, are free and open source, while others,
like MOSEK \cite{mosek}, are commercial.  These solvers can handle problems
with thousands of variables in seconds or less, and
millions of variables in minutes. Small to medium-size problems can be solved
extremely quickly using embedded solvers \cite{ecos,
stellatoOSQPOperatorSplitting2020, wang2010fast} or code generation tools
\cite{mattingley2012cvxgen, chu2013code, osqp_codegen}. For example, the
aerospace and space transportation company SpaceX uses CVXGEN
\cite{mattingley2012cvxgen} to solve convex optimization problems in real-time 
when landing the first stages of its rockets \cite{spacex}.

\paragraph{Domain-specific languages for convex optimization.}
Convex optimization problems are often specified using
domain-specific languages (DSLs) for convex optimization, such as CVXPY
\cite{diamond2016cvxpy, agrawal2018rewriting} or JuMP \cite{dunning2017jump},
which compile high-level descriptions of problems into low-level standard
forms required by solvers. The DSL then invokes a solver and retrieves a
solution on the user's behalf.  DSLs vastly reduce the engineering
effort required to get started with convex optimization, and in many cases
are fast enough to be used in production.
Using such DSLs, the convex optimization problems that we describe later
can all be implemented in just a few lines of code that very closely parallel
the mathematical specification of the problems.

\section{Constant function market makers}\label{s-cfmm}
In this section we describe how CFMMs work.
We consider a DEX with $n>1$ assets, labeled $1, \ldots, n$, that implements a CFMM.
Asset $n$ is our numeraire, the asset we use to value and assign prices to the others.

\subsection{CFMM state}\label{s-state}
\paragraph{Reserve or pool.}
The DEX has some \emph{reserves}
of available assets, given by the vector $R \in \reals_+^n$, where $R_i$ is the quantity
of asset $i$ in the reserves.
%The reserves $R$ can change over time, as 
%the \emph{transactions} described below executed.

\paragraph{Liquidity provider share weights.}
The DEX maintains a table of all the \emph{liquidity providers}, agents
who have contributed assets to the reserves.
The table includes weights
representing the fraction of the reserves each liquidity provider has a claim to.
We denote these weights as $v_1, \ldots, v_N$, where $N$ is the number of 
liquidity providers.
The weights are nonnegative and sum to one, \ie, $v \geq 0$, and $\sum_{i=1}^N v_i =1$.
The weights $v_i$ and the number of liquidity providers $N$ can change over time, 
with addition of new liquidity providers, or the deletion from the table 
of any liquidity provider whose weight is zero.

\paragraph{State of the CFMM.}
The reserves $R$ and liquidity provider weights $v$ constitute the state of the DEX.
The DEX state changes over time due to any of the three possible 
\emph{transactions}:
a \emph{trade} (or \emph{exchange}), \emph{adding liquidity},
or \emph{removing liquidity}. These transactions are described in
\S\ref{s-trade} and \S\ref{s-liquidity}.

\subsection{Proposed trade}\label{s-trade}
A \emph{proposed trade} (or \emph{proposed exchange})
is initiated by an agent or trader, 
who proposes to trade or exchange one basket of assets for another.  
A proposed trade specifies the \emph{tender basket},
with quantities given by $\Delta \in \reals_+^n$,
which is the basket of assets the trader proposes to give (or tender) to the DEX,
and the \emph{received basket},
the basket of assets the trader proposes to receive from the DEX in return,
with quantities given by $\Lambda\in \reals_+^n$.
Here $\Delta_i$ ($\Lambda_i)$ denotes the amount of asset $i$ that the trader 
proposes to tender to the DEX (receive from the DEX).
In the sequel we will refer to the vectors that give the quantities,
\ie, $\Delta$ and $\Lambda$, as the tender and receive baskets, respectively.

The proposed trade can either be rejected by the DEX, in which case
its state does not change, or accepted, in which case the basket $\Delta$ 
is transfered from the trader to the DEX, and the basket $\Lambda$ is 
transfered from the DEX to the trader.  The DEX reserves are updated as
\BEQ\label{e-update-R}
R^+ = R + \Delta - \Lambda,
\EEQ
where $R^+$ denotes the new reserves.
%When the proposed trade is accepted, we say that the trader has exchanged
%the basket $\Delta$ for the basket $\Lambda$.
A proposed trade is accepted or rejected based on a simple condition
described in \S\ref{s-function}, which always ensures that $R^+ \geq 0$.

\paragraph{Disjoint support of tender and receive baskets.}
Intuition suggests that a trade would not include an asset in both the proposed
tender and receive baskets, \ie,
we should not have $\Delta_i$ and $\Lambda_i$ both positive.
We will see later that while it is possible to include an asset 
in both baskets, it never makes sense to do so.  This means that
$\Delta$ and $\Lambda$ can be assumed to have disjoint support, \ie,
we have $\Delta_i \Lambda_i = 0$ for each $i$.
This allows us to define two disjoint sets of assets associated
with a proposed or accepted trade:
\[
\mathcal T = \{ i \mid \Delta_i >0 \}, \qquad
\mathcal R = \{ i \mid \Lambda_i >0 \}.
\]
Thus $\mathcal T$ are the indices of assets the trader proposes 
to give to the DEX, in exchange for the assets with indices in $\mathcal R$.
If $j \not\in \mathcal T \cup \mathcal R$, it means that the proposed
trade does not involve asset $j$, \ie, $\Delta_j = \Lambda_j =0$.

\paragraph{Two-asset and multi-asset trades.}
A very common type of proposed trade involves only two assets, one
that is tendered and one that is received, \ie, $|\mathcal T| =|\mathcal R|=1$.
Suppose $\mathcal T= \{i\}$ and $\mathcal R=\{j\}$, with $i\neq j$.
Then we have
$\Delta = \delta e_i$ and $\Lambda = \lambda e_j$, where $e_i$ denotes the 
$i$th unit vector, and
$\lambda \geq 0$ is the quantity of asset $j$ the trader wishes to
receive in exchange for the quantity $\delta \geq 0$ of asset $i$.
(This is referred to as exchanging asset $i$ for asset $j$.)
When a trade involves more than two assets, it is called a \emph{multi-asset}
trade. We will study two-asset and multi-asset trades in \S\ref{s-two-asset}
and \S\ref{s-multi-asset}, respectively.

\subsection{Trading function}\label{s-function}
Trade acceptance depends on both the proposed trade and the current reserves. A
proposed trade $(\Delta,\Lambda)$ is accepted only if
\BEQ\label{e-accept}
\phi(R+\gamma\Delta-\Lambda) = \phi(R),
\EEQ
where $\phi: \reals_+^n \to \reals$ is the \emph{trading function} associated with 
the CFMM, and the parameter $\gamma \in (0,1]$ introduces a \emph{trading fee}
(when $\gamma<1$). The `constant function' in the name CFMM refers to the acceptance condition \eqref{e-accept}.

We can interpret the trade acceptance condition as follows.  If $\gamma =1$, a proposed 
trade is accepted only if the quantity $\phi(R)$ does not change, \ie, $\phi(R^+)=\phi(R)$.
When $\gamma <1$ (with typical values
being very close to one), the proposed trade is accepted based on the devalued tendered 
basket $\gamma \Delta$.  The reserves, however, are updated based on the full
tendered basket $\Delta$ as in \eqref{e-update-R}.

\paragraph{Properties.}
We will assume that the trading function $\phi$ is concave, increasing, and differentiable.
Many existing CFMMs are associated with functions that satisfy the 
additional property of homogeneity, \ie, 
$\phi(\alpha R) = \alpha \phi(R)$ for $\alpha>0$.

\subsection{Trading function examples}\label{s-trading-function-examples}
We mention some trading functions that are used in existing CFMMs.

\paragraph{Linear and sum.}
The simplest trading function is linear,
\[
\phi(R) = p^TR = p_1R_1+ \cdots + p_nR_n,
\]
with $p>0$, where $p_i$ can be interpreted as the price of 
asset $i$.
The trading condition \eqref{e-accept} simplifies to 
\[
\gamma p^T \Delta = p^T \Lambda.
\]
We interpret the righthand side as the total value of received basket,
at the prices given by $p$, and the lefthand side as the 
value of the tendered basket, discounted by the factor $\gamma$.

A CFMM with $p= \ones$, \ie, all asset prices equal to one, is
called a \emph{constant sum market maker}.
The CFMM mStable, which held assets that were each pegged to the same currency,
was one of the earliest constant sum market makers.

\paragraph{Geometric mean.}
Another choice of trading function is the (weighted) geometric mean,
\[
\phi(R) = \prod_{i=1}^n R_i^{w_i},
\]
where total $w>0$ and $\ones^T w=1$. Like the linear and sum trading functions,
the geometric mean is homogeneous.

CFMMs that use the geometric mean are called \emph{constant mean market
makers}. The CFMMs Balancer \cite{balancer}, Uniswap \cite{uniswap}, and
SushiSwap \cite{sushiswap} are examples of constant mean market makers.
(Uniswap and SushiSwap use weights $w_i = 1/n$, and are sometimes called
\emph{constant product} market makers \cite{angerisAnalysisUniswapMarkets2020, angerisImprovedPriceOracles2020}.)

\paragraph{Other examples.}
Another example combines the sum and geometric mean functions,
\[
    \phi(R) = (1 - \alpha) \ones^TR + \alpha \prod_{i=1}^n R_i^{w_i},
\]
where $\alpha \in [0,1]$ is a parameter, $w \geq 0$, and $\ones^T w =1$.
This trading function yields a CFMM that
interpolates between a constant sum market (when $\alpha=0$)
and a constant geometric mean market (when $\alpha=1$).
Because it is a convex combination of the sum and geometric mean functions,
which are themselves homogeneous, the resulting function is also homogeneous.

The CFMM known as Curve
\cite{egorovStableSwapEfficientMechanism} uses the closely related trading function
\[
    \phi(R) = \ones^TR - \alpha \prod_{i=1}^n R_i^{-1},
\]
where $\alpha > 0$. Unlike the previous examples, this trading function is 
not homogeneous.

\subsection{Prices and exchange rates}\label{s-prices}
In this section we introduce the concept of asset (reported) prices,
based on a first order approximation of the trade acceptance
condition \eqref{e-accept}. These prices inform how liquidity can be
added and removed from the CFMM, as we will see in \S\ref{s-liquidity}.

%\paragraph{Num\'eraire.} It is common in practice to denote a specific asset
%to be a \emph{num\'eraire}; \ie, an asset that can be used as a common denomination
%for the value of all other assets. For the sake of convenience, we will always assume that
%asset $n$ is the numeraire, and denote all prices as being with respect to asset $n$.

\paragraph{Unscaled prices.}
We denote the gradient of the trading function as $P=\nabla \phi(R)$.
We refer to $P$, which has positive entries since $\phi$ is increasing,
as the vector of \emph{unscaled prices},
\BEQ\label{e-unscaled-prices}
P_i = \nabla \phi(R)_i = \frac{\partial \phi}{\partial R_i}(R), \quad i=1,\ldots, n.
\EEQ
To see why these numbers can be interpreted as prices, we approximate the
exchange acceptance condition \eqref{e-accept} using
its first order Taylor approximation to get
\[
0 = \phi(R+\gamma \Delta - \Lambda) - \phi(R) \approx 
\nabla \phi (R) ^T (\gamma \Delta - \Lambda)= P^T(\gamma \Delta - \Lambda),
\]
when $\gamma \Delta - \Lambda$ is small, relative to $R$.
We can express this approximation as
\BEQ\label{e-unscaled-price-interp}
\gamma \sum_{i \in \mathcal T} P_i \Delta_i
\approx \sum_{i \in \mathcal R} P_i \Lambda_i.
\EEQ
The righthand side is the value of the received basket using the unscaled prices
$P_i$.  The lefthand side is the value of the tendered basket
using the unscaled prices $P_i$, discounted by the factor $\gamma$.

\paragraph{Prices.}
The condition \eqref{e-unscaled-price-interp} is homogeneous in the prices, \ie, it is the
same condition if we scale all prices by any positive constant.
The \emph{reported prices} (or just \emph{prices}) 
of the assets are the prices relative to the price of the numeraire, which is 
asset $n$.
The prices are
\[
p_i = \frac{P_i}{P_n}, \quad i=1, \dots, n.
\]
(The price of the numeraire is always $1$.)
In general the prices depend on the reserves $R$.  (The one exception is 
with a linear trading function, in which the prices are constant.)
In terms of prices, the condition \eqref{e-unscaled-price-interp} is
\BEQ\label{e-price-interp}
\gamma \sum_{i \in \mathcal T} p_i \Delta_i
\approx \sum_{i \in \mathcal R} p_i \Lambda_i.
\EEQ

We observe for future use that the prices for two values of the reserves $R$
and $\tilde R$ are the same if and only if
\BEQ\label{e-same-price}
\nabla \phi(\tilde R) = \alpha \nabla \phi(R),
\EEQ
for some $\alpha >0$.

\paragraph{Geometric mean trading function prices.} 
For the special case
$\phi(R)=\prod_{i=1}^n R_i^{w_i}$, with $w_i >0$ and $\sum_{i=1}^nw_i=1$,
the unscaled prices are
\[
P = \nabla\phi(R) = \phi(R) (w_1R_1^{-1}, w_2R_2^{-1}, \dots, w_nR_n^{-1}),
\]
and the prices are
\begin{equation}\label{e-price-mean}
    p_i = \frac{w_iR_n}{w_nR_i},
\quad i=1, \ldots, n.
\end{equation}

\paragraph{Exchange rates.}
In a two-asset trade with $\Delta= \delta e_i$ and $\Lambda = \lambda e_j$,
\ie, we are exchanging asset $i$ for asset $j$, the \emph{exchange rate} is
\[
E_{ij} = \gamma \frac{\nabla \phi(R)_i}{\nabla \phi(R)_j} = \gamma \frac{P_i}{P_j}
= \gamma \frac{p_i}{p_j}.
\]
This is approximately how much asset $j$ you get for each unit of asset $i$, for 
a small trade.
Note that $E_{ij}E_{ji} = \gamma^2 < 1$, when $\gamma<1$, \ie, round-trip 
trades lose value.

%\paragraph{Numeraire.}
%The linearized trading acceptance condition \eqref{e-price-interp} 
%is an approximation of what baskets can be trade.  It is homogeneous
%in the prices, so the condition is the same if scale all the prices by 
%any positive number.
%We can scale the prices so that one of 
%the assets, the \emph{numeraire}, has price one.

\paragraph{These are first order approximations.}
We remind the reader that the various conditions described above
are based on a first order Taylor 
approximation of the trade acceptance condition.
A proposed trade that satisfies \eqref{e-price-interp} is
not (quite) valid; it is merely close to valid when the proposed
trade baskets are small compared to the reserves.
This is similar to the midpoint price (average of bid and ask prices)
in an order book;
you cannot trade in either direction exactly at this price.

\paragraph{Reserve value.} The value of the reserves (using the prices $p$) 
is given by
\BEQ\label{e-R-value}
V = p^T R = \frac{\nabla \phi(R)^TR}{\nabla \phi(R)_n}.
\EEQ
When $\phi$ is homogeneous we can use the identity $\nabla\phi(R)^T R =\phi(R)$
to express the reserves value as
\BEQ\label{e-R-value-homo}
V = p^T R = \frac{\phi(R)}{\nabla \phi(R)_n}.
\EEQ

\subsection{Adding and removing liquidity}\label{s-liquidity}
In this section we describe how agents called \emph{liquidity providers}
can add or remove liquidity from the reserves.
When an agent adds liquidity, she adds a basket $\Psi \in \reals_+^n$
to the reserves, resulting in the updated reserves $R^+ = R+\Psi$.
When an agent removes liquidity, she removes a basket $\Psi \in \reals_+^n$
from the reserves, resulting in the updated reserves $R^+ = R-\Psi$.
(We will see below that the condition for removing liquidity ensures that 
$R^+ \geq 0$.)
Adding or removing liquidity also updates the liquidity provider share
weights, as described below.

\paragraph{Liquidity change condition.}
Adding or removing liquidity must be done in a way that preserves the asset prices.
Using \eqref{e-same-price}, this means we must have
\BEQ\label{e-liquidity-change}
\nabla \phi(R^+) = \alpha \nabla \phi(R),
\EEQ 
for some $\alpha>0$.  (We will see later that $\alpha>1$ corresponds to removing
liquidity, and $\alpha<1$ corresponds to adding liquidity.)
This liquidity change condition is analogous to the trade exchange condition
\eqref{e-accept}.
We refer to $\Psi$ as a \emph{valid liquidity change} if this condition holds.

The liquidity change condition \eqref{e-liquidity-change} simplifies in some
cases.
For example, with a linear trading function the prices are constant, 
so any basket can be used to add
liquidity, and any basket with $\Psi \leq R$ can be removed.  (The constraint
comes from the requirement
$R^+ \geq 0$, the domain of $\phi$.)

\paragraph{Liquidity change condition for homogeneous trading function.}
Another simplification occurs when the trading function is homogeneous.
For this case we have, for any $\alpha >0$,
\[
\nabla \phi(\alpha R)= \nabla \phi(R),
\]
(by taking the gradient of $\phi(\alpha R) = \alpha \phi(R)$
with respect to $R$).
This means that $\Psi = \nu R$, for $\nu >0$, is a valid liquidity change
(provided $\nu \leq 1$ for liquidity removal).
In words: you can add or remove liquidity by adding or removing 
a basket proportional to the current reserves. 
% I find this notation rather confusing... we are taking the gradient of \phi(\alpha R) with respect
% to R, and evaluating it at R, which differs from the usual definition where \nabla\phi(\alpha R) is the
% gradient of \phi evaluated at \alpha R. Perhaps this is just me, though? XXX

\paragraph{Liquidity provider share update.}
Let $V = p^T R$ denote the value of the reserves before the liquidity change,
and $V^+=(p^+)^T R^+ =  p^T R^+$ the value after.
The change in reserve value is $V^+-V= p^T \Psi$ when adding liquidity, 
and $V^+- V= -p^T\Psi$ when removing liquidity.
Equivalently, $p^T\Psi$ is the value of the basket a liquidity provider gives,
when adding liquidity, or receives when removing liquidity.
The fractional change in reserve value is $(V^+-V)/V^+$.

When liquidity provider $j$ adds or removes liquidity, all the share weights 
are adjusted pro-rata based on the change of value of the reserves, 
which is the value of the basket she adds or removes.
The weights are adjusted to
\BEQ \label{e-share-update}
v_i^+ = \begin{cases}
    v_i V/V^+ + (V^+ - V)/V^+ & i= j\\
    v_iV/V^+ & i \ne j.
\end{cases}
\EEQ
Thus the weight of liquidity provider $j$ is increased (decreased) 
by the fractional change in reserve value when she adds (removes) liquidity.
These new weights are also nonnegative and sum to one.

When $\phi$ is homogeneous and we add liquidity with the basket $\Psi = \nu R$,
with $\nu >0$, we have $V_+ = (1+\nu)p^TR$, so 
\[
V/V^+ = 1/(1+\nu), \qquad (V^+-V)/V^+ =  \nu / (1+\nu).
\]
The weight updates for adding liquidity $\Psi= \nu R$ are then
\[
v_i^+ = \begin{cases}
(v_i+\nu)/(1+\nu) & i=j \\
v_i/(1+\nu) &  i\neq j.
\end{cases}
\]
For removing liquidity with the basket $\Psi = \nu R$,
we replace $\nu$ with $-\nu$ in the formulas above, along with
the constraint $\nu \le v_j$.

\subsection{Agents interacting with CFMMs}\label{s-agents-interacting}
Agents seeking to trade or add or remove liquidity make proposals.
These proposals are accepted or not, depending on the acceptance conditions 
given above. A proposal can be rejected if another agent's proposed action is 
accepted (processed) before their proposed action, thus changing $R$ and invalidating
the acceptance condition.

\paragraph{Slippage thresholds.} \label{p-slippage}
One practical and common approach to mitigating this problem during trading is to
allow agents to set a \emph{slippage threshold} on the received basket. This slippage threshold,
represented as some
percentage $0 \le \eta \le 1$, is simply a parameter that specifies how much slippage
the agent is willing to tolerate without their trade failing. In this case,
the agent presents some trade $(\Delta, \Lambda)$ along with a threshold $\eta$,
and the contract accepts the trade if there is some number $\alpha$ satisfying
$\eta \le \alpha$ such that the trade $(\Delta, \alpha\Lambda)$ can be accepted.
In other words, the agent allows the contract to devalue the output basket by
at most a factor of $\eta$. If no such value of
$\alpha$ exists, the trade fails.

\paragraph{Maximal liquidity amounts.} While setting slippage thresholds can help
with reducing the risk of trades failing, another possible failure mode can occur
during the addition of liquidity. 
A simple solution to this problem is that the liquidity
provider specifies some basket $\Psi$ to the CFMM contract, and the contract accepts
the largest possible basket $\Psi^-$ such that $\Psi^- \le \Psi$, returning
the remaining amount, $\Psi - \Psi^-$, to the liquidity provider. In other words,
$\Psi$ can be seen as the maximal amount of liquidity a user is willing to provide.

\section{Properties}\label{s-prop}
In this section we present some basic properties of CFMMs.

\subsection{Properties of trades}
\paragraph{Non-uniqueness.} 
If we replace the trading function $\phi$ with $\tilde \phi = h\circ \phi$, where
$h$ is concave, increasing, and differentiable, we obtain another concave increasing 
differentiable function.
The associated CFMM has the same trade acceptance condition, 
the same prices, the same liquidity change condition, and the same 
liquidity provider share updates as the original CFMM.

\paragraph{Maximum valid receive basket.}
Any valid trade satisfies
$\phi(R+\gamma \Delta-\Lambda) = \phi(R)$, 
so in particular $R+\gamma \Delta -\Lambda \geq 0$.
Since we assume $\Delta$ and $\Lambda$ have non-overlapping support, it follows that
\[
\Lambda \le R.
\]
A valid trade cannot ask to receive more than is in the reserves.

\paragraph{Non-overlapping support for valid tender and receive baskets.}
Here we show why a valid proposed trade with $\Delta_k>0$ and $\Lambda_k>0$ for some $k$
does not make sense when $\gamma<1$, justifying our assumption that this never happens.
Let $(\tilde \Delta,\tilde \Lambda)$ be a proposed trade which coincides with
$(\Delta,\Lambda)$ except in the $k$th components, which we set to
\[
\tilde \Delta_k = \Delta_k - \tau/\gamma, \qquad
\tilde \Lambda_k = \Lambda_k - \tau,
\]
where $\tau = \min\{\gamma \Delta_k,\Lambda_k\} >0$.
Evidently $\tilde \Delta \geq 0$, $\tilde \Lambda \geq 0$, and
\[
R+\gamma \Delta - \Lambda = R+\gamma \tilde \Delta - \tilde \Lambda,
\]
so the proposed trade $(\tilde \Delta,\tilde \Lambda)$ is also valid.
If the trader proposes this trade instead of $(\Delta,\Lambda)$, the net change
in her assets is 
\[
\tilde \Lambda- \tilde \Delta = \Lambda -\Delta + 
\left(\frac{1}{\gamma}-1\right) \tau e_k.
\]
The last vector on the right is zero in all entries except $k$,
and positive in that entry.  Thus the valid proposed trade $(\tilde \Delta,\tilde \Lambda)$
has the same net effect as the trade $(\Delta,\Lambda)$, except that the trader
ends up with a positive amount more of the $k$th asset.  Assuming the 
$k$th asset has value, we would always prefer this.

\paragraph{Trades increase the function value.} 
For an accepted nonzero trade, we have
\[
\phi(R^+) =\phi(R+\Delta-\Lambda) >
\phi(R+\gamma\Delta- \Lambda) = \phi(R),
\]
since $\phi$ is increasing and $R+\Delta-\Lambda \ge R+\gamma\Delta- \Lambda$, 
with at least one entry being strictly greater, whenever $\gamma < 1$.

We can derive a stronger inequality using concavity of $\phi$, which implies that
\[
\phi(R+\gamma \Delta-\Lambda) \leq \phi(R+\Delta-\Lambda) +
(\gamma-1) \nabla \phi(R+\Delta-\Lambda)^T \Delta.
\]
This can be re-arranged as
\[
\phi(R^+) \geq \phi(R) + (1-\gamma) (P^+)^T \Delta,
\]
where $P^+ = \nabla \phi(R^+)$ are the unscaled prices at the reserves $R^+$.
This tells us the function value increases at least by $(1-\gamma)$ times
the value of tendered basket at the unscaled prices.

%If $\phi$ is homogeneous, we divide this inequality by 
%$\nabla \phi(R^+)_n$ abd use \eqref{e-R-value-homo} to express it as
%\[
%(p^+)^TR^+ \geq \frac{\phi(R)}{P^+_n} +  (1-\gamma)(p^+)^T \Delta.
%\]

\paragraph{Trading cost is positive.}
Suppose $(\Delta, \Lambda)$ is a valid trade.
The net change in the trader's holdings is $\Lambda-\Delta$.  We can interpret
$\delta = p^T (\Delta-\Lambda)$ as the decrease in value of the trader's holdings
due to the proposed trade, evaluated at the current prices.  We can interpret
$\delta$ as a trading cost, evaluated at the pre-trade prices,
and now show it is positive.

Since $\phi$ is concave, we have
\[
\phi(R+\gamma \Delta -\Lambda) \leq \phi(R) + \nabla \phi(R)^T(\gamma \Delta - \Lambda).
\]
Using $\phi(R+\gamma \Delta -\Lambda) = \phi(R)$, this implies
\[
0 \leq \nabla \phi(R)^T (\gamma \Delta - \Lambda) = P^T(\gamma \Delta-\Lambda).
%\sum_{i=1}^n P_i (\gamma \Delta_i - \Lambda_i).
\]
From this we obtain
\[
P^T(\Delta-\Lambda) = P^T(\gamma \Delta-\Lambda) + (1-\gamma) P^T \Delta 
\geq (1-\gamma) P^T \Delta .
\]
Dividing by $P_n$ gives
\[
\delta \geq (1-\gamma) p^T \Delta.
\]
Thus the trading cost is always at least a factor $(1-\gamma)$ of $p^T\Delta$,
the total value of the tendered basket.

The trading cost $\delta$ is also the \emph{increase} in the total reserve value,
at the current prices.
So we can say that each trade increases the total reserve value, at the current prices,
by at least $(1-\gamma)$ times the value of the tendered basket.

%\paragraph{No splitting of trades.} An important consequence of the increasing constant property is that a trader will never be incentivized to split
%up trades. More specifically, we will show that if $(\Delta, \Lambda)$ and $(\Delta + \Delta^+, \Lambda + \Lambda^+)$ are both valid trades
%for the current reserve values, then, executing $(\Delta, \Lambda)$ and attempting to trade $(\Delta^+, \Lambda^+)$ after will always be underbid; \ie, it can
%only be accepted if the trader either adds more input or requests less output.
%
%This follows nearly immediately from the definition. Since $(\Delta + \Delta^+, \Lambda + \Lambda^+)$ is valid, we have:
%\[
%\phi(R + \gamma(\Delta + \Delta^+) - (\Lambda + \Lambda^+)) = \phi(R).
%\]
%On the other hand, we have that:
%\[
%\phi(R + \gamma(\Delta + \Delta^+) - (\Lambda + \Lambda^+)) < \phi(R + \Delta - \Lambda + \gamma \Delta^+ - \Lambda^+)
%\]
%
%XXX: Tried for some amount of time and couldn't really simplify it. Can include it later if you want to chat, but it does need a few definitions set up

\subsection{Properties of liquidity changes}\label{s-liquidity-properties}

\paragraph{Liquidity change condition interpretation.}
One natural interpretation of the liquidity change condition~\eqref{e-liquidity-change}
is in terms of a simple optimization problem.
We seek a basket $\Psi$ that maximizes the
post-change trading function value subject to a given total value of the basket
at the current prices,
\begin{equation}\label{e-liquidity-problem}
\begin{aligned}
    & \text{maximize} && \phi(R^+)\\
    & \text{subject to} && p^T (R^+-R) \le M.
\end{aligned}
\end{equation}
Here the optimization variable is $R^+\in \reals_+^n$, and
$M$ is the desired value of the basket $\Psi$ at the current prices,
for adding liquidity, or its negative, for removing liquidity.
The optimality conditions for this convex optimization problem are
\[
p^T(R^+-R)\le M, \qquad \nabla \phi(R^+) - \nu p =0,
\]
where $\nu \ge 0$ is a Lagrange multiplier. Using $p = \nabla \phi(R)/\nabla \phi(R)_n$, the second condition is 
\[
\nabla \phi(R^+) = \frac{\nu}{\nabla \phi(R)_n} \nabla \phi(R),
\]
which is \eqref{e-liquidity-change} with $\alpha = \nu/ \nabla \phi(R)_n$.
We can easily recover the trading basket $\Psi$ from $R^+$ since $\Psi = R^+ - R$.

\paragraph{Liquidity provision problem.} When the trading function is homogeneous,
it is easy to understand what baskets can be used to add or remove liquidity:
they must be proportional to the current reserves.  In other cases, it can be 
difficult to find an $R^+$ that satisfies \eqref{e-liquidity-change}.
In the general case, however,
the convex optimization problem~\eqref{e-liquidity-problem} can be solved 
to find the basket $\Psi$ that gives a valid liquidity change, with $M$
denoting the total value of the added basket (when $M>0$) or 
removed basket (when $M<0$).

%\paragraph{Maximal liquidity amounts.}
%As described in~\S\ref{s-agents-interacting}, because the reserves of the CFMM are
%constantly changing due to the actions of other agents, an agent attempting to add a specific basket $\Psi$ of liquidity,
%which can be found by solving~\eqref{e-liquidity-problem}, might fail. In this case, the CFMM
%can instead solve the following problem:
%\[
%\begin{aligned}
%    & \text{maximize} && \phi(R^+)\\
%    & \text{subject to} && p^T (R^+-R) \le p^T\Psi\\
%    &&& R^+ - R \le \Psi,
%\end{aligned}
%\]
%with variable $R^+$. The data in this problem are the prices $p$, the current reserves
%$R^+$ and the maximal basket of assets $\Psi$, \ie, the maximum amount of each asset which the
%liquidity provider wishes to contribute.

\paragraph{Liquidity change and the gradient scale factor $\alpha$.}
Suppose that we add or remove liquidity.
Since $\phi$ is concave \eqref{e-monotone} tells us that
\[
(\nabla \phi(R^+)-\nabla \phi(R))^T (R^+-R) \leq 0.
\]
Using $\nabla \phi(R^+) = \alpha \nabla \phi(R)$, this becomes
\[
(\alpha-1) \nabla \phi(R)^T (R^+-R) \leq 0.
\]
We have $\nabla \phi(R)>0$.  
If we add liquidity, we have 
$R^+ - R \geq 0$ and $R^+-R \neq 0$, so $\nabla \phi(R)^T (R^+-R) >0$.
From the inequality above we conclude that $\alpha<1$.
If we remove liquidity, a similar arguments tells us that $\alpha>1$.

\section{Two-asset trades}\label{s-two-asset}
Two-asset trades, sometimes called \emph{swaps}, are some of the most common
types of trades performed on DEXs. In this section, we show a number of interesting properties
of trades in this common special case.

\subsection{Exchange functions}

Suppose we exchange asset $i$ for asset $j$, so $\Delta = \delta e_i$ and 
$\Lambda = \lambda e_j$, with $\delta\geq 0$, $\lambda \geq 0$.
The trade acceptance condition \eqref{e-accept} is 
\BEQ \label{e-accept-2}
\phi(R+\gamma \delta e_i - \lambda e_j)=\phi(R).
\EEQ
The lefthand side is increasing in $\delta$ and decreasing in $\lambda$,
so for each value of $\delta$ there is at most one valid value of $\lambda$,
and for each value of $\lambda$, there is at most one valid value of
$\delta$.
In other words, the relation \eqref{e-accept-2} 
between $\lambda$ and $\gamma$ defines a one-to-one function.
This means that two-asset trades are characterized by a single parameter,
either $\delta$ (how much is tendered) or $\lambda$ (how much is received).

\paragraph{Forward exchange function.}
Define $F: \reals_+ \to \reals$,
where $F(\delta)$ is the unique $\lambda$ that satisfies \eqref{e-accept-2}.
The function $F$ is called the \emph{forward exchange function},
since $F(\delta)$ is how much of asset $j$ you get if you exchange
$\delta$ of asset $i$. The forward exchange function $F$ is increasing since $\phi$ is 
componentwise increasing and nonnegative since $F(0) = 0$. 
We will now show that the function $F$ is concave.

\paragraph{Concavity.} Using the implicit function theorem on~\eqref{e-accept-2}
with $\lambda = F(\delta)$, we obtain
\BEQ\label{e-exchange-scalar}
F'(\delta) = \gamma\frac{\nabla\phi(R')_i}{\nabla\phi(R')_j},
\EEQ
where we use $R' = R + \gamma\delta e_i - F(\delta) e_j$ to simplify notation.
%If $\delta = 0$ then $F'(0) = \gamma P_i/P_j = E_{ij}$ which is simply the exchange
%rate between asset $i$ and $j$.
To show that $F$ is concave, we will show that, for any nonnegative trade amounts $\delta, \delta' \ge 0$, the
function $F$ satisfies
\BEQ\label{e-concavity-scalar}
F(\delta') \le F'(\delta)(\delta' - \delta) + F(\delta),
\EEQ
which establishes that $F$ is concave.

We write $R'' = R + \gamma \delta' e_i - F(\delta')e_j$, and note that 
$\phi(R) = \phi(R') = \phi(R'')$ from the definition of $F$.
Since $\phi$ is concave it satisfies
\[
\phi(R'') \le \nabla \phi(R')^T(R'' - R') + \phi(R'),
\]
so $\nabla \phi(R')^T(R'' - R') \ge 0$. 
Using the definitions of $R''$ and $R'$, we have
\[
0 \le \gamma(\delta' - \delta)\nabla \phi(R')_i - (F(\delta') - F(\delta))\nabla\phi(R')_j.
\]
Dividing by $\nabla\phi(R')_j$ and using~\eqref{e-exchange-scalar},
we obtain \eqref{e-concavity-scalar}.

\paragraph{Reverse exchange function.}
Define $G: \reals_+ \to \reals \cup \{\infty\}$,
where $G(\lambda)$ is the unique $\delta$ that satisfies \eqref{e-accept-2},
or $G(\lambda)=\infty$ is there is no such $\delta$.
The function $F$ is called the \emph{reverse exchange function},
since $F(\lambda)$ is how much of asset $i$ you must exchange, to
receive $\lambda$ of asset $j$. In a similar way to the forward trade function,
the reverse exchange function is nonnegative 
and increasing, but this function is convex rather than concave.
(This follows from a nearly identical proof.) 

\paragraph{Forward and reverse exchange functions are inverses.}
The forward and reverse exchange functions are inverses of each other,
\ie, they satisfy
\[
G(F(\delta)) = \delta, \qquad F(G(\lambda)) = \lambda,
\]
when both functions are finite.

\paragraph{Analogous functions for a limit order book market.}
There are analogous functions in a market that uses a limit order
book.  They are piecewise linear, where the slopes are the different prices
of each order, while the distance between the kink points is equal to the size 
of each order.  The associated functions have the same properties,
\ie, they are increasing, inverses of each other, $F$ is concave, and $G$ is 
convex.

\paragraph{Evaluating $F$ and $G$.}
In some important special cases, we can express the 
functions $F$ and $G$ in a closed form.
For example, when the trading function is the sum function, they are
\[
F(\delta) = \min\{\gamma\delta, R_j\}, \qquad G(\lambda) = \begin{cases}
    \lambda/\gamma & \lambda/\gamma \le R_j\\
    +\infty & \text{otherwise}.
\end{cases}
\]
When the trading function is the geometric mean, the functions are
\[
F(\delta) = R_j \left(1-\frac{R_i^{w_i/w_j}}{(R_i + \gamma \delta)^{w_i/w_j}}\right), \qquad
G(\lambda) = \frac{R_i}{\gamma}\left(\frac{R_j^{w_j/w_i}}{(R_j - \lambda)^{w_j/w_i}} - 1\right),
\]
whenever $\lambda < R_j$, and $G(\lambda) = \infty$ otherwise.

On the other hand, when the forward and reverse trading functions $F$ and $G$
cannot be expressed analytically, we can use several methods
to evaluate them numerically~\cite[\S9]{press1992numerical}.
To evaluate $F(\delta)$, we fix $\delta$ and solve for $\lambda$ 
in \eqref{e-accept-2}.   The lefthand side is a decreasing function of
$\lambda$, so we can use simple bisection to solve this nonlinear equation.
Newton's method can be used to achieve higher accuracy with fewer steps.
Exploiting the concavity of $\phi$, it can be shown an undamped
Newton iteration always converges to the solution.
With superscripts denoting iteration, this is
\[
\lambda^{k+1} = \lambda^k + \frac{\phi(R+\gamma \delta e_i - \lambda^k e_j) - \phi(R)}
{\nabla \phi(R+\gamma \delta e_i - \lambda^k e_j)_j},
\]
with starting point based on the exchange rate,
\[
\lambda^0 = \delta E_{ij} = \delta\frac{\gamma p_i}{p_j}.
\]
(It can be shown that the convergence is monotone decreasing.)
We note that one of the largest CFMMs, Curve, uses a trading function that
is not homogeneous and uses this method
in production \cite{egorovStableSwapEfficientMechanism}.

\paragraph{Slope at zero.}
Using~\eqref{e-exchange-scalar}, we see that $F'(0^+) = E_{ij}$,
\ie, the one-sided derivative at $0$ is exactly the exchange rate for 
assets $i$ and $j$.
Since $F$ is concave, we have
\BEQ\label{e-price-upper-bound}
F(\delta) \leq  F'(0^+)\delta = E_{ij} \delta.
\EEQ
This tells us that the amount of asset $j$ you will receive for trading $\delta$ of 
asset $i$ is no more than the amount predicted by the exchange rate.

The one-sided derivative of the reverse exchange function $G$
at $0$ is $G'(0^+) = E_{ji}$.  The analog of 
the inequality \eqref{e-price-upper-bound} is
\BEQ\label{e-price-lower-bound}
G(\lambda) \geq  G'(0^+)\lambda = \gamma^{-2} E_{ji} \lambda,
\EEQ
which states that the amount of asset $i$ you need to tender to receive 
an amount of asset $j$ is at least the amount predicted by the exchange rate.

\begin{figure}
    \centering
    \includegraphics[width=.98\textwidth]{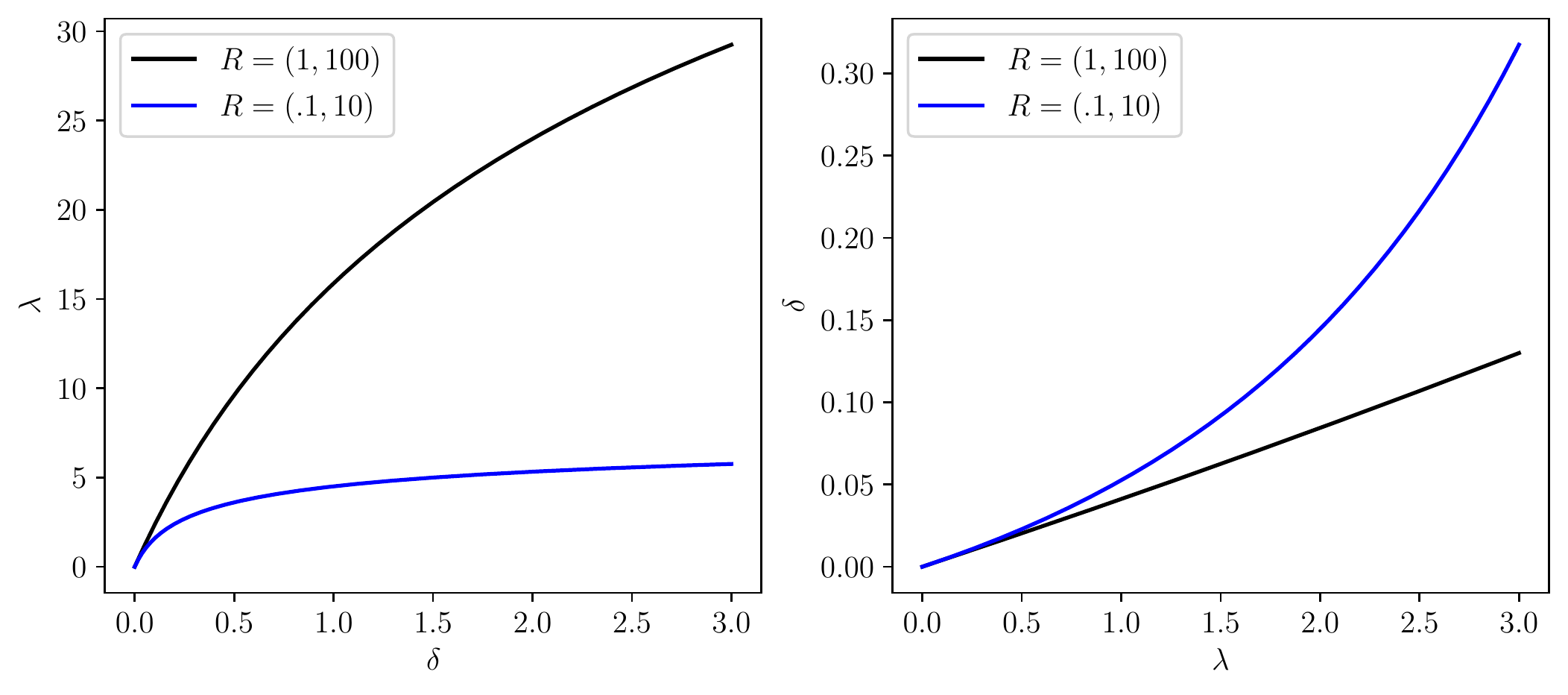}
    \caption{\emph{Left.} Forward exchange functions for two 
values of the reserves. \emph{Right.} Reverse exchange functions
for the same two values of the reserves.}
    \label{f-forward}
\end{figure}

\paragraph{Examples.}
Figure~\ref{f-forward} shows the forward and reverse exchange functions
for a constant geometric mean market with two assets and weights $w_1 = .2$
and $w_2 = .8$, and $\gamma =0.997$.
We show the functions for two values of the reserves:
$R=(1,100)$ and $R=(0.1,10)$.
The exchange rate is
the same for both values of the reserves and equal to 
$E_{12} = \gamma w_1R_2/w_2R_1 = 25$.

\subsection{Exchanging multiples of two baskets}
Here we discuss a simple generalization of two-asset trade, 
in which we tender and receive a multiple of fixed baskets.
Thus, we have $\Delta = \delta \tilde \Delta$ and $\Lambda = \lambda \tilde \Lambda$,
where $\lambda\geq 0$ and $\delta \geq 0$ scale the fixed baskets $\tilde \Delta$
and $\tilde \Lambda$.
When $\tilde \Delta = e_i$ and $\tilde \Lambda= e_j$, this reduces to the 
two-asset trade discussed above.

The same analysis holds in this case as in the simple two-asset trade.
We can introduce the forward and reverse functions $F$ and $G$, which are inverses
of each other. They are increasing, $F$ is concave, $G$ is convex,
and they satisfy $F(0)=G(0)=0$.
We have the inequality
\[
F(\delta) \leq E \delta, 
\]
where $E$ is the exchange rate for exchanging the basket $\tilde \Delta$ for the 
basket $\tilde \Lambda$, given by 
\[
E = \gamma \frac{\nabla \phi(R)^T \tilde \Delta}{\nabla \phi(R)^T \tilde \Lambda}.
\]
There is also an inequality analogous to \eqref{e-price-lower-bound}, using
this definition of the exchange rate. We mention two specific 
important examples in what follows.

\paragraph{Liquidating assets.}
Let $\Delta\in \reals_+^n$ denote a basket of assets we wish to liquidate,
\ie, exchange for the numeraire.  We can assume that $\Delta_n=0$.
We then find the $\alpha>0$ for which $(\Delta,\alpha e_n)$ is a valid
trade, \ie,
\begin{equation}\label{e-liquidation-value}
\phi(R+\gamma \Delta - \alpha e_n) = \phi(R).
\end{equation}
We can interpret $\alpha$ as the \emph{liquidation value} of the basket $\Delta$      .
We can also show that the liquidation value is at most as large as the discounted
value of the basket; \ie, $\alpha \le \gamma p^T\Delta$.

To see this, apply~\eqref{e-concavity} to the left hand side of~\eqref{e-liquidation-value},
which gives, after cancelling $\phi(R)$ on both sides,
\[
\nabla\phi(R)^T(\gamma \Delta - \alpha e_n) \ge 0.
\]
Rearranging, we find:
\[
\alpha \le \frac{\gamma\nabla\phi(R)^T\Delta}{\nabla\phi(R)_n} = \gamma p^T\Delta.
\]

%XXX: add liquidated value, show that liquidation value always decreases after trades

\paragraph{Purchasing a basket.}
Let $\Lambda \in \reals_+^n$ denote a basket we wish to purchase using 
the numeraire.  We find $\alpha >0$ for which
$(\alpha e_n,\Lambda)$ is a valid trade, \ie,
\[
\phi(R+\gamma \alpha e_n - \Lambda) = \phi(R).
\]
We interpret $\alpha$ as the \emph{purchase cost} of the basket $\Lambda$.
It can be shown that $\alpha \geq (1/\gamma) p^T \Lambda$, \ie, the purchase cost
is at least a factor $1/\gamma$ more than the value of the 
basket, at the current prices. This follows from a nearly identical
argument to that of the liquidation value.

\section{Multi-asset trades}\label{s-multi-asset}
We have seen that two-asset trades are easy to understand; we
choose the amount we wish to tender (or receive), and we can then find
the amount we will receive (or tender).
Multi-asset trade are more complex, because even for a fixed receive basket
$\Lambda$, there are many tender baskets that are valid, and we face the question
of which one should we use.  The same is true when we fix the tendered
basket $\Delta$: there are many baskets $\Lambda$ we could receive, 
and we need to choose one.
More generally, we have the question of how to choose the   
proposed trade $(\Delta,\Lambda)$.
In the two-asset case, the choice is parametrized by a scalar, 
either $\delta$ or $\lambda$.  In the multi-asset case, there are more degrees 
of freedom.

\begin{figure}
    \centering
    \includegraphics[width=.8\textwidth]{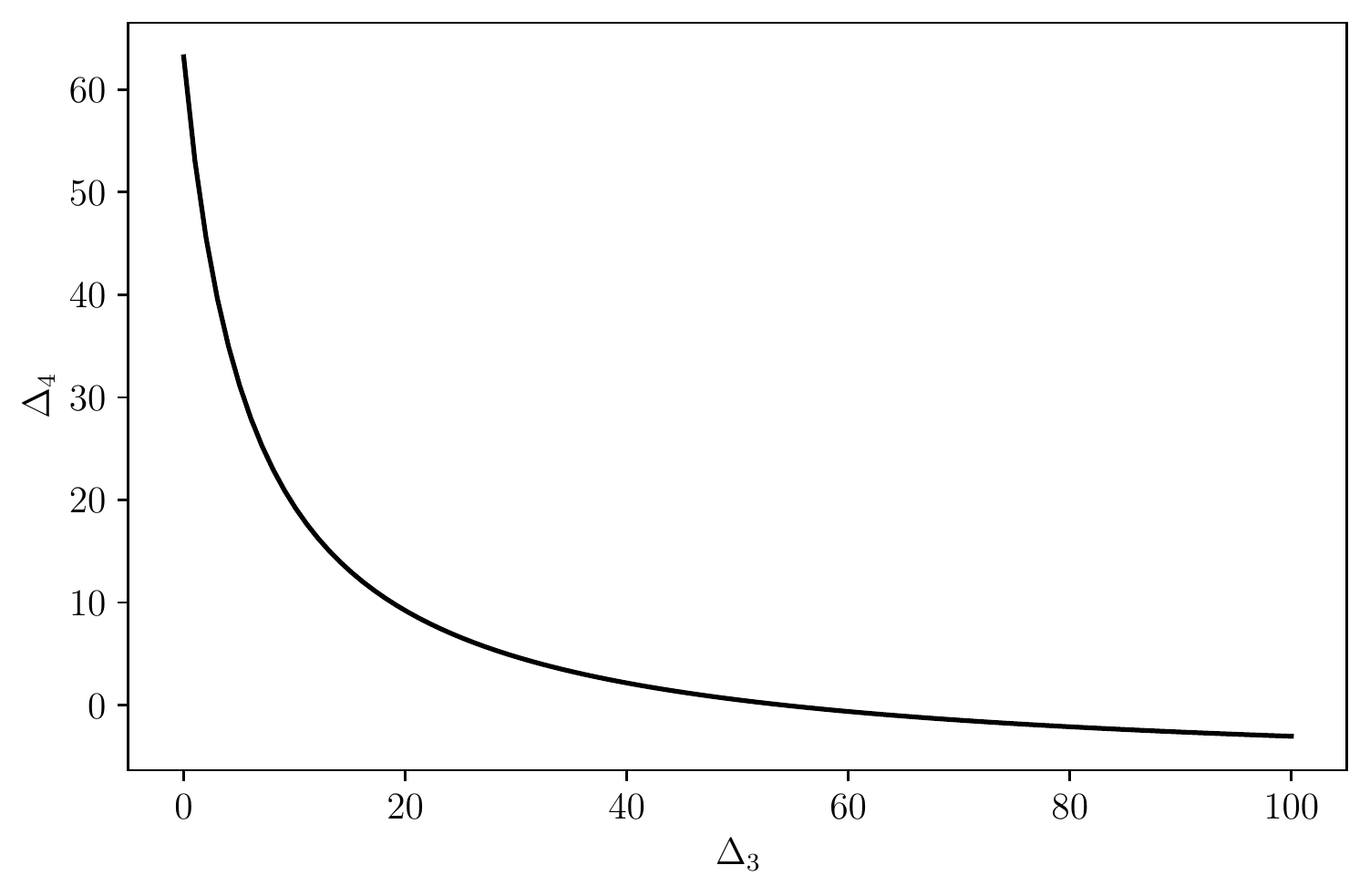}
    \caption{Valid tendered baskets $(\Delta_3, \Delta_4)$ for the 
        received basket $\Lambda =(2,4,0,0)$.}
    \label{f-tendered}
\end{figure}

\paragraph{Example.}
We consider an example with $n=4$, geometric mean trading function with 
weights $w_i = 1/4$ and fee $\gamma = .997$, with reserves $R = (4, 5, 6, 7)$.
We fix the received basket to be $\Lambda = (2,4,0,0)$.
There are many valid tendered baskets, which
are shown in figure~\ref{f-tendered}. 
The plot shows valid values of $(\Delta_3, \Delta_4)$,
since the first two components of $\Delta$ are zero.

\subsection{The general trade choice problem}
We formulate the problem of choosing $(\Delta,\Lambda)$ as an optimization problem.
The net change in holdings of the trader is $\Lambda - \Delta$.
The trader judges a net change in holdings using a utility function $U:\reals^n \to
\reals \cup\{-\infty\}$, where
she prefers $(\Delta,\Lambda)$ to $(\tilde \Delta, \tilde \Lambda)$ if
$U(\Lambda - \Delta)> U(\tilde \Lambda - \tilde \Delta)$.
The value $-\infty$ is used to indicate that a change in holdings is
unacceptable.
We will assume that $U$ is increasing and concave.
(Increasing means that the trader would always prefer to have a larger 
net change than a smaller one, which comes from our assumption
that all assets have value.)

To choose a valid trade that maximizes utility, we solve the problem
\BEQ\label{e-prob}
\begin{aligned}
& \text{maximize} && U(\Lambda - \Delta) \\
& \text{subject to} && \phi(R+\gamma \Delta - \Lambda)  = \phi(R), \quad 
\Delta \geq 0, \quad \Lambda \geq 0,
\end{aligned}
\EEQ
with variables $\Delta$ and $\Lambda$. 
Unfortunately the constraint $\phi(R+\gamma \Delta - \Lambda) = \phi(R)$ is not 
convex (unless the trading function is linear),
so this problem is not in general convex.

Instead we will solve its convex relaxation, where we
change the equality constraint to an inequality to obtain the convex problem
\BEQ\label{e-prob-cvx}
\begin{aligned}
& \text{maximize} && U(\Lambda - \Delta) \\
& \text{subject to} && \phi(R+\gamma \Delta - \Lambda) \geq \phi(R), \quad 
\Delta \geq 0, \quad \Lambda \geq 0,
\end{aligned}
\EEQ
which is readily solved.
It is easy to show that any solution of \eqref{e-prob-cvx} satisfies
$\phi(R+\gamma \Delta - \Lambda) = \phi(R)$, and so is also a solution of
the problem \eqref{e-prob}.
(If a solution satisfies
$\phi(R+\gamma \Delta - \Lambda) > \phi(R)$, we can decrease $\Delta$ or increase 
$\Lambda$ a bit, so as to remain feasible and increase the objective, 
a contradiction.)

Thus we can (globally and efficiently) 
solve the non-convex problem \eqref{e-prob} by solving
the convex problem \eqref{e-prob-cvx}.

\paragraph{No-trade condition.}
Assuming $U(0) > - \infty$, the solution to the problem \eqref{e-prob-cvx} can 
be $\Delta = \Lambda = 0$, which means that
trading does not increase the trader's utility, \ie, the trader should not 
propose any trade.
We can give simple conditions under which this happens for the case
when $U$ is differentiable.
They are
\BEQ\label{e-no-trade}
\gamma p \leq \alpha\nabla U(0) \leq p,
\EEQ
for some $\alpha>0$. We can interpret the set of prices $p$ for which this is true, \ie,
\[
K = \{p \in \reals^n_+ \mid \gamma p \le \alpha \nabla U(0) \le p ~ \text{for some} ~ \alpha > 0\},
\]
as the \emph{no-trade cone} for the utility function $U$. 
(It is easy to see that $K$ is a convex polyhedral cone.)

We interpret $\nabla U(0)$ as the vector of marginal utilities to the trader, 
and $p$ as the prices of the assets in the CFMM.  
For $\gamma =1$, the condition says that we do not trade when
the marginal utility is a positive multiple of the current asset prices; 
if this does not hold, then the solution of the trading problem
\eqref{e-prob-cvx} is nonzero, \ie, the trader should trade to increase 
her utility.
When $\gamma<1$, the trader will not trade when the prices are in $K$.

To derive condition~\eqref{e-no-trade}, we first derive the optimality conditions for
the problem \eqref{e-prob-cvx}.  We introduce the Lagrangian
\[
L(\Delta, \Lambda, \lambda, \omega, \kappa) =
U(\Lambda -\Delta) + \lambda(\phi(R+\gamma \Delta-\Lambda)-\phi(R)) +
\omega^T\Delta + \kappa^T\Lambda,
\]
where $\lambda \in \reals_+$, $\omega\in \reals_+^n$, and 
$\kappa \in \reals_+^n$ are dual variables or Lagrange multipliers for
the constraints.
The optimality conditions for \eqref{e-prob-cvx} are feasibility, along with
\[
\nabla_\Delta L = 0, \qquad \nabla _\Lambda L = 0.
\]
The choice $\Delta=0$, $\Lambda =0$ is feasible, and satisfies this
condition if
\[
\nabla_\Delta L(0,0,\lambda,\omega,\kappa) = 0, \qquad 
\nabla_\Lambda L(0,0,\lambda,\omega,\kappa) = 0.
\]
These are
\[
- \nabla U(0) + \lambda \gamma \nabla \phi(R) + \omega =0, \qquad
\nabla U(0) - \lambda \nabla \phi(R) + \kappa=0,
\]
which we can write as
\[
\nabla U(0) \geq \lambda \gamma \nabla \phi(R), \qquad
\nabla U(0) \leq \lambda \nabla \phi(R).
\]
Dividing these by $\lambda P_n$, we obtain \eqref{e-no-trade},
with $\alpha =  1/(\lambda P_n)$.

\subsection{Special cases}
\paragraph{Linear utility.}
When $U(z) = \pi^T z$, with $\pi \geq 0$, we can interpret $\pi$ as the trader's private
prices of the assets, \ie, the prices she values the assets at.
From \eqref{e-no-trade} we see that the trader
will not trade if her private asset prices satisfy
\BEQ\label{e-no-arb}
\gamma p \le \alpha\pi \le p
\EEQ
for some $\alpha > 0$.

In the special case where $\pi$ satisfies
\[
(\pi_2, \dots, \pi_n) = \lambda (p_2, \dots, p_n),
\]
for $\lambda \ge 0$, \ie, $\pi$ is collinear with $p$ except in the first entry,
\eqref{e-no-arb} is satisfied if and only if
\[
\lambda \gamma p_1 \le \pi_1 \le \lambda \gamma^{-1}p_1.
\]
If $\lambda = 1$, then this simplifies to the condition
\[
\gamma p_1 \le \pi_1 \le \gamma^{-1}p_1.
\]
(This will arise in an example we present below.)

\paragraph{Markowitz trading.}
Suppose the trader models the return $r\in \reals^n$ on the assets over
some period of time as a random vector with mean 
$\Expect r = \mu \in \reals^n$ and 
covariance matrix $\Expect (r-\mu)(r-\mu)^T =\Sigma \in \reals^{n\times n}$.
If the trader holds a portfolio of assets $z \in \reals^n_+$, the return
is $r^T z$;
the expected portfolio return is $\mu^Tz$ and the variance of the portfolio return
is $z^T\Sigma z$.
In Markowitz trading, the trader maximizes the 
risk-adjusted return, defined as $\mu^T z - \kappa z^T\Sigma z$,
where $\kappa>0$ is the 
\emph{risk-aversion parameter}~\cite{markowitzPortfolioSelection1952, boyd2017multi}.
This leads to the Markowitz trading problem
\BEQ\label{e-markowitz}
\begin{aligned}
& \mbox{maximize} && \mu^T z - \kappa z^T \Sigma z\\
&\mbox{subject to} && z = z^\text{curr} - \Delta + \Lambda \\
&&&\phi(R+\gamma \Delta - \Lambda) \geq \phi(R)\\
&&& \Delta \geq 0, \quad \Lambda \geq 0,
\end{aligned}
\EEQ
with variables $z$, $\Delta$, $\Lambda$, where $z^\text{curr}$ is the trader's 
current holdings of assets.
This is the general problem \eqref{e-prob-cvx} with concave utility function
\[
U(Z) = \mu^T(z^\text{curr} + Z) - \kappa
(z^\text{curr} + Z)^T \Sigma (z^\text{curr} + Z).
\]

A well-known limitation of the Markowitz quadratic utility function $U$, \ie, the
risk-adjusted return, is that
it is not increasing for all $Z$, which implies that the trading function
relaxation need not be tight.  However, for any sensible choice 
of the parameters $\mu$ and $\Sigma$, it is increasing for the values of
$Z$ found by solving the Markowitz problem \eqref{e-markowitz},
and the relaxation is tight.
As a practical matter, if a solution of \eqref{e-markowitz} does not satisfy 
the trading constraint, then the parameters are inappropriate.

\paragraph{Expected utility trading.}
Here the trader models the returns $r\in \reals^m$ on the assets over
some time interval as random, with some known distribution.
The trader seeks to maximize the
expected utility of the portfolio return, using a concave increasing 
utility function $\psi: \reals \to \reals$ to introduce risk aversion.
(Thus we use the term utility function to refer to both the trading utility function
$U: \reals_+^n \to \reals$ and the portfolio 
return utility function $\psi:\reals \to \reals$,
but the context should make it clear which is meant.)
This leads to the problem
\BEQ\label{e-max-util-trading}
\begin{aligned}
& \mbox{maximize} && \Expect \psi(r^Tz)\\
&\mbox{subject to} && z = z^\text{curr} - \Delta + \Lambda \\
&&&\phi(R+\gamma \Delta - \Lambda) \geq \phi(R)\\
&&& \Delta \geq 0, \quad \Lambda \geq 0,
\end{aligned}
\EEQ
where the expectation is over $r$.  This is the general problem
\eqref{e-prob-cvx}, with utility
\[
U(Z) = \Expect \psi (r^T (z^\text{curr}+Z)),
\]
which is concave and increasing.

This problem can be solved using several methods. 
One simple approach is to replace the expectation
with an empirical or sample average over some Monte Carlo 
samples of $r$, which leads to an approximate solution of
\eqref{e-max-util-trading}.   The problem can also be solved 
using standard methods for convex stochastic optimization, 
such as projected stochastic gradient methods.

\subsection{Numerical examples}
In this section we give two numerical examples.

\paragraph{Linear utility.}
Our first example involves a CFMM with 6 assets,
geometric mean trading function with equal weights $w_i = 1/6$, 
and trading fee parameter $\gamma = .9$. (We intentionally use an unrealistically
small value of $\gamma$ so the no-trade condition is more evident.)
We take reserves
\[
R = (1, 3, 2, 5, 7, 6).
\]
The corresponding prices are given by~\eqref{e-price-mean},
\[
p = (R_6/R_1, R_6/R_2, \ldots, 1) = (6, 2, 3, 6/5, 6/7, 1).
\]
We consider linear utility, with the trader's private prices given by
\[
\pi = (tp_1, p_2, \dots, p_n),
\]
where $t$ is a parameter that we vary over the interval $t \in [1/2,2]$. 
For $t=1$, we have $\pi=p$, \ie, the CFMM prices and the trader's private
prices are the same (and not surprisingly, the trader does not trade).
As we vary $t$, we vary the trader's private price for asset $1$ by up to
a factor of two from the CFMM price.

\begin{figure}
    \centering
    \includegraphics[width=.99\textwidth]{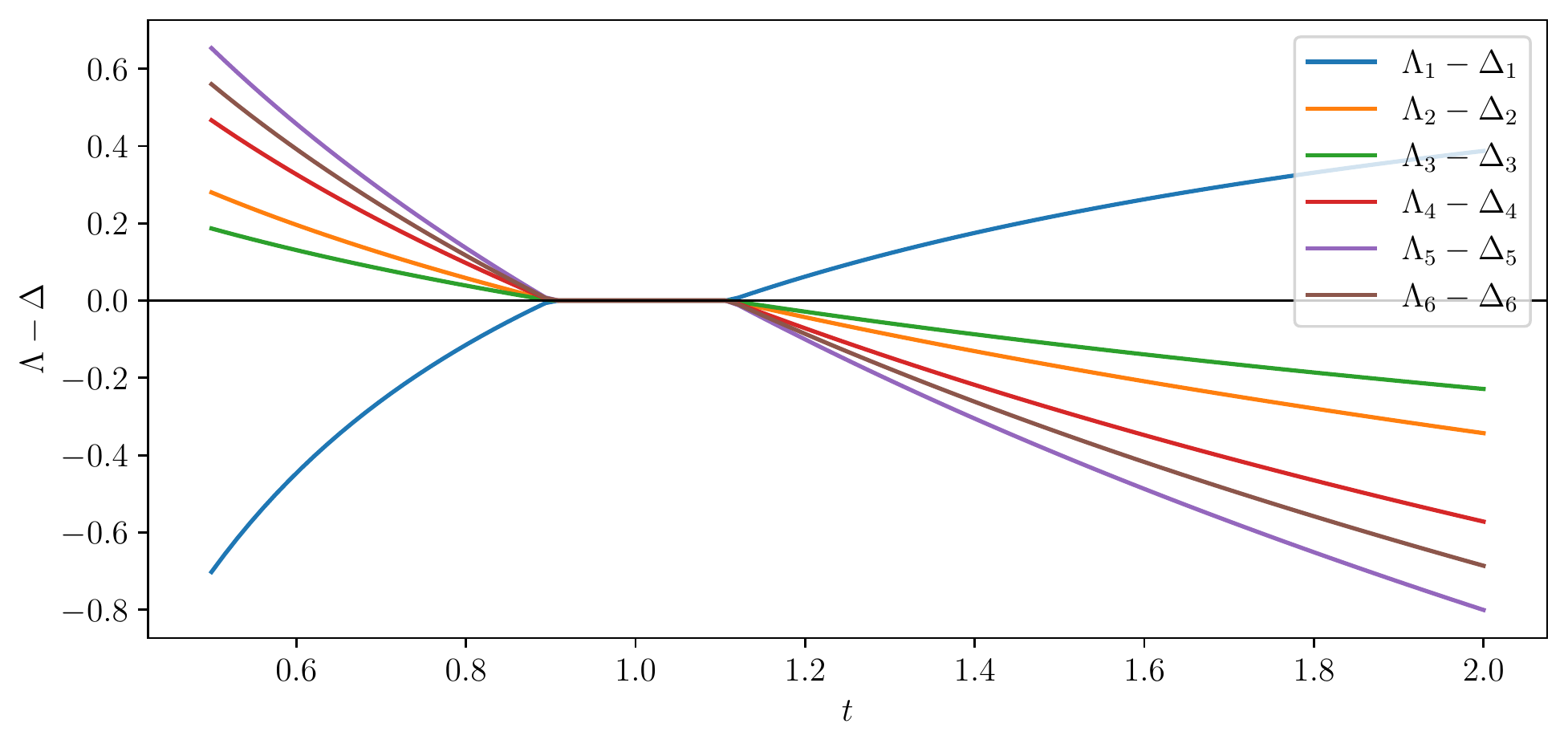}
    \caption{Solutions $\Lambda - \Delta$ for the linear utility 
maximization problem, as the private price for asset $1$ is varied 
by the factor $t$ from the CFMM price. The blue curve shows asset $1$.}
    \label{f-linear-utility}
\end{figure}

The family of optimal trades are shown in
figure~\ref{f-linear-utility}, as a function of the parameter $t$.
We plot $\Lambda - \Delta$ versus $t$, which shows assets in the tender
basket as negative and the received basket as positive.
The blue curve shows asset $1$, which we tender when $t$ is small,
and receive when $t$ is large.
The no-trade region is clearly seen as the interval $t\in [0.9,1.1]$.

\begin{figure}
    \centering
    \includegraphics[width=.99\textwidth]{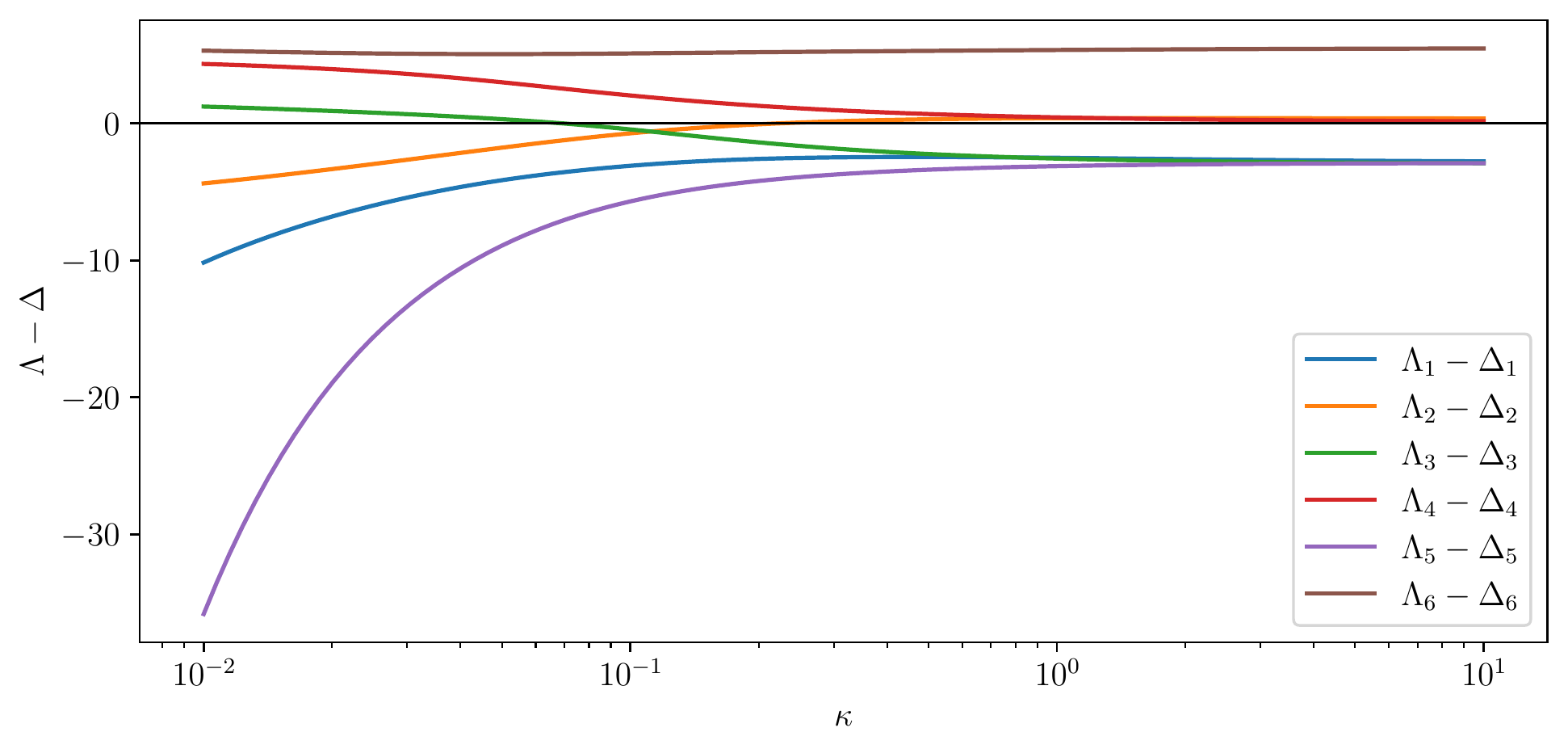}
    \caption{Solutions $\Lambda - \Delta$ for instances of an example Markowitz trading problem
    as the risk-aversion parameter $\kappa$ is varied.}
    \label{f-markowitz}
\end{figure}

\paragraph{Markowitz trading.} 
Our second example uses nearly the same CFMM and reserves as the previous example,
but with a more realistic trading fee parameter
$\gamma = .997$. (This is a common choice of trading fee for many CFMMs.)
%In this example, we take
%\[
%R = (10, 30, 20, 50, 70, 60),
%\]
%which are just the reserves in the previous example, scaled by a factor of 10, and take
We solve the Markowitz trading problem \eqref{e-markowitz}, with current holdings
\[
z^\text{curr} = (2.5, 1, .5, 2.5, 3, 1),
\]
mean return
\[
\mu = (-.01, .01, .03, .05, -.02, .02),
\]
and covariance $\Sigma = V^TV/100$,
where the entries of $V\in \reals^{6\times 6}$ 
are drawn from the standard normal distribution.
We solve the optimal trading problem for values of the risk aversion parameter 
$\kappa$ varying between $10^{-2}$ and $10^{1}$.
(For all of these values, the trading constraint is tight.)
These optimal trades are shown in figure~\ref{f-markowitz}.
It is interesting to note that depending on the risk aversion, we either
tender or receive assets 2 and 3.

The CVXPY code for the Markowitz optimal trading problem is given below.
In this snippet we assume that \verb|mu|, \verb|sigma|, \verb|gamma|, 
\verb|kappa|, \verb|R|, and \verb|z_curr| have been previously defined.
Note that the code closely follows the mathematical description of the problem
given in \eqref{e-markowitz}.

\begin{lstlisting}[float=h, language=mypython, caption=Markowitz trading CVXPY code.]
import cvxpy as cp

delta = cp.Variable(6)
lam = cp.Variable(6)

z = z_curr - delta + lam
R_new = R + gamma*delta - lam

objective = cp.Maximize(z.T @ mu - kappa*cp.quad_form(z, sigma))
constraints = [
    cp.geo_mean(R_new) >= cp.geo_mean(R),
    delta >= 0,
    lam >= 0
]

problem = cp.Problem(objective, constraints)
problem.solve()
\end{lstlisting}

\section{Conclusion}

We have provided a general description of CFMMs, outlining how users can
interact with a CFMM through trading or adding and removing liquidity. 
We observe that many of the properties of CFMMs follow from concavity of
the trading function.
In the simple case where two assets are traded or exchanged,
it suffices to specify the amount we wish
to receive (or tender), which determines the amount we
tender (receive), by simply evaluating a convex (concave) function. 
Multi-asset trades are more complex, since the set of valid trades is 
multi-dimensional, \ie, multiple tender or received baskets are possible.
We formulate the problem of choosing from among these possible valid trades
as a convex optimization problem, which can be globally and efficiently solved.

\clearpage

\section*{Acknowledgements}
The authors would like to acknowledge Shane Barratt for useful discussions.
Guillermo Angeris is supported by the National Science Foundation Graduate Research Fellowship under Grant No.\ DGE-1656518.
Akshay Agrawal is supported by a Stanford Graduate Fellowship.

%\section{Liquidity changes}\label{s-liquidity-changes}
%XXX it's a one dimensional problem, specified by $\alpha$r or equivalently by
%the post transaction weights $w^+$.
%
%XXX plot shows the basket you have to tender versus $\alpha$, and also the 
%post-transcation weight $w^+$.
%
%XXX can also parametrize the basket versus $w^+$.
%
%So this problem is like a two asset trade; it depends on only one 
%scalar parameter.
%
%The problem is more interesting when combined with trades either before or after 
%the liquidity change.   These can be solved by solving a one parameter problem
%that characterizes the liquidity change, and then solving a convex problem
%for the trade.  (The cost function can include not only the posterior holdings,
%but also the posterior weight.)

\bibliographystyle{alpha}
\bibliography{references.bib}

\newcommand{\etalchar}[1]{$^{#1}$}
\begin{thebibliography}{WCDW21}

\bibitem[aav21]{aave}
Aave.
\newblock https://aave.com, 2021.

\bibitem[ABN{\etalchar{+}}21]{agrawal2021allocation}
Akshay Agrawal, Stephen Boyd, Deepak Narayanan, Fiodar Kazhamiaka, and Matei
  Zaharia.
\newblock Allocation of fungible resources via a fast, scalable price discovery
  method.
\newblock {\em arXiv preprint arXiv:2104.00282}, 2021.

\bibitem[AC20]{angerisImprovedPriceOracles2020}
Guillermo Angeris and Tarun Chitra.
\newblock Improved {{price oracles}}: {{Constant function market makers}}.
\newblock In {\em Proceedings of the 2nd {{ACM Conference}} on {{Advances}} in
  {{Financial Technologies}}}, pages 80--91, {New York NY USA}, October 2020.
  {ACM}.

\bibitem[AEC20]{angeris2020does}
Guillermo Angeris, Alex Evans, and Tarun Chitra.
\newblock When does the tail wag the dog? {C}urvature and market making.
\newblock {\em arXiv preprint arXiv:2012.08040}, 2020.

\bibitem[AEC21a]{angeris2021note}
Guillermo Angeris, Alex Evans, and Tarun Chitra.
\newblock A note on privacy in constant function market makers.
\newblock {\em arXiv preprint arXiv:2103.01193}, 2021.

\bibitem[AEC21b]{angeris2021replicating}
Guillermo Angeris, Alex Evans, and Tarun Chitra.
\newblock Replicating market makers.
\newblock {\em arXiv preprint arXiv:2103.14769}, 2021.

\bibitem[AI21]{aoyagiLiquidityImplicationsConstant}
Jun Aoyagi and Yuki Ito.
\newblock Liquidity implications of constant product market makers.
\newblock {\em Available at SSRN 3808755}, 2021.

\bibitem[AKC{\etalchar{+}}20]{angerisAnalysisUniswapMarkets2020}
Guillermo Angeris, Hsien-Tang Kao, Rei Chiang, Charlie Noyes, and Tarun Chitra.
\newblock An {{analysis}} of {{Uniswap}} markets.
\newblock {\em Cryptoeconomic Systems}, November 2020.

\bibitem[Aoy20]{aoyagiLiquidityProvisionAutomated}
Jun Aoyagi.
\newblock Liquidity provision by automated market makers.
\newblock {\em Available at SSRN 3674178}, 2020.

\bibitem[ApS19]{mosek}
MOSEK ApS.
\newblock {{MOSEK Optimizer API}} for {{Python}} 9.1.5.
\newblock https://docs.mosek.com/9.1/pythonapi/index.html, 2019.

\bibitem[AVB21]{angerisHeuristicMethodsPerformance2021}
Guillermo Angeris, Jelena Vu{\v c}kovi{\'c}, and Stephen Boyd.
\newblock Heuristic methods and performance bounds for photonic design.
\newblock {\em Optics Express}, 29(2):2827, January 2021.

\bibitem[AVDB18]{agrawal2018rewriting}
Akshay Agrawal, Robin Verschueren, Steven Diamond, and Stephen Boyd.
\newblock A rewriting system for convex optimization problems.
\newblock {\em Journal of Control and Decision}, 5(1):42--60, 2018.

\bibitem[AZS{\etalchar{+}}21]{adams2021uniswap}
Hayden Adams, Noah Zinsmeister, Moody Salem, River Keefer, and Dan Robinson.
\newblock Uniswap v3 core.
\newblock Technical report, 2021.

\bibitem[BBD{\etalchar{+}}17]{boyd2017multi}
Stephen Boyd, Enzo Busseti, Steven Diamond, Ronald Kahn, Kwangmoo Koh, Peter
  Nystrup, and Jan Speth.
\newblock Multi-period trading via convex optimization.
\newblock {\em Foundations and Trends in Optimization}, 3(1):1--76, 2017.

\bibitem[BKPH05]{boyd2005digital}
Stephen Boyd, Seung-Jean Kim, Dinesh Patil, and Mark Horowitz.
\newblock Digital circuit optimization via geometric programming.
\newblock {\em Operations Research}, 53(6), 2005.

\bibitem[Bla16]{spacex}
Lars Blackmore.
\newblock Autonomous precision landing of space rockets.
\newblock {\em The BRIDGE}, 26(4), 2016.

\bibitem[BPC{\etalchar{+}}11]{boyd2011distributed}
Stephen Boyd, Neal Parikh, Eric Chu, Borja Peleato, and Jonathan Eckstein.
\newblock Distributed optimization and statistical learning via the alternating
  direction method of multipliers.
\newblock {\em Foundations and Trends{\textregistered} in Machine learning},
  3(1):1--122, 2011.

\bibitem[BSM{\etalchar{+}}17]{osqp_codegen}
Goran Banjac, Bartolomeo Stellato, Nicholas Moehle, Paul Goulart, Alberto
  Bemporad, and Stephen Boyd.
\newblock Embedded code generation using the {OSQP} solver.
\newblock In {\em {IEEE} Conference on Decision and Control}, 2017.

\bibitem[But13]{buterin2013ethereum}
Vitalik Buterin.
\newblock Ethereum: A next-generation smart contract and decentralized
  application platform, 2013.

\bibitem[But17]{buterin2017path}
Vitalik Buterin.
\newblock On path independence.
\newblock https://vitalik.ca/general/2017/06/22/marketmakers.html, 2017.

\bibitem[BV04]{boyd2004convex}
Stephen Boyd and Lieven Vandenberghe.
\newblock {\em Convex Optimization}.
\newblock {Cambridge University Press}, {Cambridge, UK ; New York}, 2004.

\bibitem[CAEK21]{chitra2021note}
Tarun Chitra, Guillermo Angeris, Alex Evans, and Hsien-Tang Kao.
\newblock A note on borrowing constant function market maker shares.
\newblock 2021.

\bibitem[CFL{\etalchar{+}}08]{chen2008complexity}
Yiling Chen, Lance Fortnow, Nicolas Lambert, David Pennock, and Jennifer
  Wortman.
\newblock Complexity of combinatorial market makers.
\newblock In {\em Proceedings of the 9th ACM Conference on Electronic
  Commerce}, pages 190--199, 2008.

\bibitem[com21]{compound}
Compound.
\newblock https://compound.finance, 2021.

\bibitem[CPDB13]{chu2013code}
Eric Chu, Neal Parikh, Alexander Domahidi, and Stephen Boyd.
\newblock Code generation for embedded second-order cone programming.
\newblock In {\em European Control Conference}, pages 1547--1552. IEEE, 2013.

\bibitem[CT06]{cornuejols2006}
Gerard Cornuejols and Reha T\"{u}t\"{u}nc\"{u}.
\newblock {\em Optimization Methods in Finance}.
\newblock Cambridge University Press, 2006.

\bibitem[DB16]{diamond2016cvxpy}
Steven Diamond and Stephen Boyd.
\newblock {CVXPY}: A {P}ython-embedded modeling language for convex
  optimization.
\newblock {\em Journal of Machine Learning Research}, 17(83):1--5, 2016.

\bibitem[DCB13]{ecos}
Alexander Domahidi, Eric Chu, and Stephen Boyd.
\newblock {{ECOS}}: {{An SOCP}} solver for embedded systems.
\newblock In {\em 2013 {{European Control Conference}} ({{ECC}})}, pages
  3071--3076, {Zurich}, July 2013. {IEEE}.

\bibitem[DHL17]{dunning2017jump}
Iain Dunning, Joey Huchette, and Miles Lubin.
\newblock {{JuMP}}: A modeling language for mathematical optimization.
\newblock {\em SIAM review}, 59(2):295--320, 2017.

\bibitem[dyd21]{dydx}
dydx.
\newblock https://dydx.exchange, 2021.

\bibitem[EAC21]{evans2021optimal}
Alex Evans, Guillermo Angeris, and Tarun Chitra.
\newblock Optimal fees for geometric mean market makers.
\newblock {\em arXiv preprint arXiv:2104.00446}, 2021.

\bibitem[Ego19]{egorovStableSwapEfficientMechanism}
Michael Egorov.
\newblock {{StableSwap}} - efficient mechanism for {{Stablecoin}} liquidity.
\newblock page~6, 2019.

\bibitem[Eva20]{evans2020liquidity}
Alex Evans.
\newblock Liquidity provider returns in geometric mean markets.
\newblock {\em arXiv preprint arXiv:2006.08806}, 2020.

\bibitem[FHT01]{friedman2001elements}
Jerome Friedman, Trevor Hastie, and Robert Tibshirani.
\newblock {\em The Elements of Statistical Learning}, volume~1.
\newblock Springer Series in Statistics, 2001.

\bibitem[GCG19]{garstkaCOSMOConicOperator2019}
Michael Garstka, Mark Cannon, and Paul Goulart.
\newblock {{COSMO}}: {{A}} conic operator splitting method for large convex
  problems.
\newblock In {\em 2019 18th {{European Control Conference}} ({{ECC}})}, pages
  1951--1956, {Naples, Italy}, June 2019. {IEEE}.

\bibitem[Han03]{hanson2003combinatorial}
Robin Hanson.
\newblock Combinatorial information market design.
\newblock {\em Information Systems Frontiers}, 5(1):107--119, 2003.

\bibitem[HBL01]{hershenson2001opamp}
Maria Hershenson, Stephen Boyd, and Thomas Lee.
\newblock Optimal design of a {CMOS} op-amp via geometric programming.
\newblock {\em IEEE Transactions on Computer-aided design of integrated
  circuits and systems}, 20(1):1--21, 2001.

\bibitem[LB14]{lipp2014minimum}
Thomas Lipp and Stephen Boyd.
\newblock Minimum-time speed optimisation over a fixed path.
\newblock {\em International Journal of Control}, 87(6):1297--1311, 2014.

\bibitem[Lu17]{lu2017building}
Alan Lu.
\newblock Building a decentralized exchange in {E}thereum.
\newblock
  https://blog.gnosis.pm/building-a-decentralized-exchange-in-ethereum-eea4e7452d6e,
  2017.

\bibitem[Mar52]{markowitzPortfolioSelection1952}
Harry Markowitz.
\newblock Portfolio {{selection}}.
\newblock {\em The Journal of Finance}, 7(1):77--91, 1952.

\bibitem[MB12]{mattingley2012cvxgen}
Jacob Mattingley and Stephen Boyd.
\newblock {CVXGEN}: A code generator for embedded convex optimization.
\newblock {\em Optimization and Engineering}, 13(1):1--27, 2012.

\bibitem[MBBW19]{moehle2018dynamicnot}
Nicholas Moehle, Enzo Busseti, Stephen Boyd, and Matt Wytock.
\newblock Dynamic energy management.
\newblock {\em arXiv preprint arXiv:1903.06230}, 2019.

\bibitem[MM19]{balancer}
Fernando Martinelli and Nikolai Mushegian.
\newblock Balancer: {{A}} non-custodial portfolio manager, liquidity provider,
  and price sensor.
\newblock 2019.

\bibitem[Nak08]{nakamoto2008bitcoin}
Satoshi Nakamoto.
\newblock Bitcoin: A peer-to-peer electronic cash system, 2008.

\bibitem[OCPB16]{odonoghueConicOptimizationOperator2016}
Brendan O'Donoghue, Eric Chu, Neal Parikh, and Stephen Boyd.
\newblock Conic {{optimization}} via {{operator splitting}} and {{homogeneous
  self}}-{{dual embedding}}.
\newblock {\em Journal of Optimization Theory and Applications},
  169(3):1042--1068, June 2016.

\bibitem[PTFV92]{press1992numerical}
William Press, Saul Teukolsky, Brian Flannery, and William Vetterling.
\newblock {\em Numerical Recipes: The Art of Scientific Computing}.
\newblock Cambridge University Press, 1992.

\bibitem[RB16]{ryuPrimerMonotoneOperator2016}
Ernest Ryu and Stephen Boyd.
\newblock A primer on monotone operator methods.
\newblock {\em Applied Computational Math}, 2016.

\bibitem[SB08]{stewart2008model}
Gregory Stewart and Francesco Borrelli.
\newblock A predictive control framework for industrial turbodiesel engine
  control.
\newblock In {\em {IEEE} Conference on Decision and Control ({CDC})}, pages
  5704--5711, 2008.

\bibitem[SBG{\etalchar{+}}20]{stellatoOSQPOperatorSplitting2020}
Bartolomeo Stellato, Goran Banjac, Paul Goulart, Alberto Bemporad, and Stephen
  Boyd.
\newblock {{OSQP}}: An operator splitting solver for quadratic programs.
\newblock {\em Mathematical Programming Computation}, February 2020.

\bibitem[Sus20]{sushiswap}
Sushi.
\newblock The {SushiSwap} project, 2020.

\bibitem[Sza95]{szabo1995smart}
Nick Szabo.
\newblock Smart contracts.
\newblock {\em Extropy: Journal of Transhumanist Thought}, 16, 1995.

\bibitem[TW20]{tassy2020growth}
Martin Tassy and David White.
\newblock Growth rate of a liquidity provider’s wealth in $xy= c$ automated
  market makers, 2020.

\bibitem[uma21]{uma}
{UMA} project.
\newblock https://umaproject.org, 2021.

\bibitem[WB10]{wang2010fast}
Yang Wang and Stephen Boyd.
\newblock Fast evaluation of quadratic control-{L}yapunov policy.
\newblock {\em IEEE Transactions on Control Systems Technology},
  19(4):939--946, 2010.

\bibitem[WCDW21]{wang2021cyclic}
Ye~Wang, Yan Chen, Shuiguang Deng, and Roger Wattenhofer.
\newblock Cyclic arbitrage in decentralized exchange markets.
\newblock {\em Available at SSRN 3834535}, 2021.

\bibitem[Win69]{winkler1969scoring}
Robert Winkler.
\newblock Scoring rules and the evaluation of probability assessors.
\newblock {\em Journal of the American Statistical Association},
  64(327):1073--1078, 1969.

\bibitem[Woo14]{wood2014ethereum}
Gavin Wood.
\newblock Ethereum: A secure decentralised generalised transaction ledger,
  2014.

\bibitem[Woo16]{wood2016polkadot}
Gavin Wood.
\newblock Polkadot: {V}ision for a heterogeneous multi-chain framework, 2016.

\bibitem[Yak18]{yakovenko2018solana}
Anatoly Yakovenko.
\newblock Solana: A new architecture for a high performance blockchain, 2018.

\bibitem[ZCP18]{uniswap}
Yi~Zhang, Xiaohong Chen, and Daejun Park.
\newblock Formal specification of constant product ($xy=k$) market maker model
  and implementation.
\newblock 2018.

\end{thebibliography}

\end{document}